\documentclass[paper=a4, abstract=true, twocolumn=false, twoside=semi, fontsize=12pt]{scrartcl}

\usepackage[utf8]{inputenc}
\usepackage[T1]{fontenc}
\usepackage[english]{babel}

\usepackage{amsmath}
\usepackage{amsthm}
\usepackage{mathtools}
\usepackage{nicefrac}

\usepackage{graphicx}
\usepackage[dvipsnames,table]{xcolor}
\usepackage{tikz}

\usetikzlibrary{fit,calc,shadings}

\usepackage{booktabs}
\usepackage{tabularx}
\usepackage{bigdelim}
\usepackage{multirow}

\newcolumntype{g}{>{\columncolor{gray!25}}c} 

\usepackage{enumitem}
\usepackage{ragged2e}
\usepackage{relsize}
\usepackage{microtype}
\usepackage{scrlayer-scrpage}
\usepackage{xspace}
\usepackage[breaklinks=true,colorlinks=true,bookmarks=false,linkcolor=BrickRed,citecolor=Green,hypertexnames=false]{hyperref}
\usepackage{bm} 
\usepackage{etoolbox}
\usepackage[all]{nowidow}

\usepackage[backend=biber,style=trad-plain]{biblatex}
\DeclareFieldFormat*{title}{#1}
\DeclareFieldFormat*{titlecase}{#1}
\usepackage[ruled,vlined,nofillcomment,linesnumbered]{algorithm2e}

\SetCommentSty{MyCommentStyle}
\SetProgSty{}
\SetAlCapHSkip{0pt}
\usepackage[capitalize,noabbrev]{cleveref}

\makeatletter
\newcommand\Paragraph{\@startsection{paragraph}{4}{\z@}%
                                    {0pt plus .2em}%
                                    {-.4em plus -1.0em}%
                                    {\normalfont\normalsize\bfseries}}
\newcommand\ParagraphRun{\@startsection{paragraph}{4}{\z@}%
                                       {0pt plus .2em}%
                                       {-\fontdimen2\font plus -\fontdimen3\font minus -\fontdimen4\font}%
                                       {\normalfont\normalsize\bfseries}}

\DeclareSectionCommand[%
  style=section,%
  level=4,%
  indent=\z@,%
  beforeskip=3.25ex \@plus1ex \@minus.2ex,%
  afterskip=-\the\fontdimen2\font plus -\the\fontdimen3\font minus -\the\fontdimen4\font,%
  tocstyle=subsection,%
  tocindentfollows=subsubsection,%
  tocindent=7.0em,%
  tocnumwidth=4.1em%
]{paragraphRunIn}
\makeatother

\setkomafont{caption}{\small\sffamily}
\setkomafont{captionlabel}{\sffamily\bfseries}
\addtokomafont{disposition}{\rmfamily}
\setcapindent{0pt}
\deffootnote{.6em}{2em}{%
	\makebox[.6em][l]{\textsuperscript{\thefootnotemark}}%
}

\usepackage{newpxtext}
\usepackage{newpxmath}

\newtheorem{theorem}{Theorem}
\newtheorem{corollary}{Corollary}
\newtheorem{remark}{Remark}
\newtheorem{proposition}{Proposition}
\newtheorem{example}{Example}
\newtheorem{definition}{Definition}
\newtheorem{lemma}{Lemma}

\makeatletter
\appto{\endtheorem}{\@endpetrue}
\appto{\endcorollary}{\@endpetrue}
\appto{\endremark}{\@endpetrue}
\appto{\endproposition}{\@endpetrue}
\appto{\endexample}{\@endpetrue}
\appto{\enddefinition}{\@endpetrue}
\appto{\endlemma}{\@endpetrue}
\makeatother

\DeclareMathOperator*\argmax{arg\,max}
\DeclareMathOperator*\Sign{Sign}

\newcommand\SC{\mathcal{C}}
\newcommand\SH{\mathcal{H}}

\newcommand\SV{\mathcal{V}}
\newcommand\SE{\mathcal{E}}
\newcommand\SG{\mathcal{G}}

\newcommand\SP{\mathcal{P}}

\newcommand\BR{\mathbb{R}}
\newcommand\BN{\mathbb{N}}

\newcommand\va{\mathbf{a}}
\newcommand\vb{\mathbf{b}}
\newcommand\vx{\mathbf{x}}
\newcommand\vy{\mathbf{y}}
\newcommand\vc{\mathbf{c}}

\newcommand\vu{\mathbf{u}}
\newcommand\vOne{\mathbf{1}}
\newcommand\vZero{\mathbf{0}}
\newcommand\vlambda{{\bm\lambda}}
\newcommand\valpha{{\bm\alpha}}

\newcommand\tauBatch{\tau_{\mathsf{batch}}}
\newcommand\tauDrop{\tau_{\mathsf{drop}}}
\newcommand\tauFeas{\tau_{\mathsf{feas}}}
\newcommand\tauGap{\tau_{\mathsf{gap}}}
\newcommand\tauStab{\tau_{\mathsf{stab}}}

\newcommand\eg{e.g.\xspace}
\newcommand\ie{i.e.\xspace}
\newcommand\cf{c.f.\xspace}
\newcommand\wrt{wrt.\xspace}
\newcommand\WLOG{w.l.o.g.\xspace}

\cohead{A Bregman-Sinkhorn Algorithm for the Maximum Weight Independent Set Problem}
\cehead{
  Stefan Haller,
  Bogdan Savchynskyy
}
\ofoot*{}
\ifoot*{}

\title{A Bregman-Sinkhorn Algorithm for the Maximum Weight Independent Set Problem}

\author{
  Stefan Haller,
  Bogdan Savchynskyy\\[.8ex]%
  \large
  Heidelberg University
 }

\date{}

\begin{document}

\maketitle

\begin{abstract}
  \setlength\parindent{15pt}%
  We propose a scalable approximate algorithm for the NP-hard maximum-weight independent set problem, based on dual coordinate descent applied to a smoothed clique-cover LP relaxation.
  Our method, a variant of the Bregman/Sinkhorn algorithm, employs entropy smoothing with a novel duality-gap-based smoothing scheduling strategy that empirically outperforms standard feasibility scheduling for the relaxed problem.
  Our new projection to the primal feasible set enables accurate duality gap estimation.
  Combined with a basic primal heuristic leveraging reduced costs, our approach yields high-quality integer solutions.
  On real-world datasets, it efficiently finds high-quality approximate solutions for graphs with up to 882,000~nodes and 344~million edges within seconds.
\end{abstract}

\section{Introduction}
\label{sec:intro}

Let $\SG=(\SV,\SE)$ be an undirected graph with a \emph{node set} $\SV$ and an \emph{edge set} $\SE\subseteq \binom{\SV}{2}$.
An \emph{independent} or \emph{stable set} is a subset of mutually non-adjacent vertices $\SV'\subseteq\SV$, \ie, for any $i,j\in\SV'$ it holds $\{i,j\}\notin\SE$.
Given the costs $c_i$, $i\in\SV$, the \emph{maximum weight independent set (MWIS)} problem consists in finding an independent set with the maximum total cost.

The problem has applications in computer vision~\cite{brendel2010segmentation}, \cite[Ch.~7]{prakash2022fully}, vehicle routing~\cite{dong2021new}, transmission scheduling in wireless networks~\cite{sanghavi2007message} and is one of the classical combinatorial optimization problems~\cite{garey1979computers}.
It is NP-hard since its unweighted variant is NP-complete, as shown in \cite{karp2010reducibility} for the vertex cover problem which is polynomially reducible to the maximum independent set problem.
Therefore, various relaxations and approximation algorithms addressing this problem have been proposed.

We consider a family of \emph{clique cover relaxations} for this problem and propose a scalable dual algorithm.
We build upon it with an efficient primal heuristic that profits from the reduced costs computed by the dual algorithm.
In total, this leads to an efficient and well-scalable MWIS method.

\paragraph{Notation.}
Vectors are \textbf{bold} and scalars not: $\vx$ is a vector and $x$ is a scalar.
Scalar operations and relations are applied to vectors component-wise: $\vx\ge \vZero$ means non-negativity of all coordinates of $\vx$ and $\log\vx$ is a vector of coordinate logarithms.
$\BR^n_{\ge 0} \coloneq \{\vx\in\BR^n \mid \vx\ge0\}$, $[n] \coloneq \{1,2,\dots,n\}$.
$\llbracket A \rrbracket$ are Iverson brackets, \ie, the expression is equal to $1$ if $A$ holds, otherwise 0.

\paragraph{Notation introduced later in the paper.}
$\vx$~--~primal variables, $\vc$~--~cost vector~(\cref{sec:overview});
$\vlambda$~--~dual variables, $\vc^{\vlambda}$~--~vector of reduced costs~(\cref{sec:dual-of-clique-relaxation});
$\bar K_j=K_j\cup\{x_{n+j}\}$, $j \in [m]$, are cliques with and without slack variables $x_{n+j}$~(\cref{sec:overview});
$J_i=\{j\in[m] \mid \bar K_j\ni i\}$~--~the set of cliques containing variable $x_i$~(\cref{sec:dual-of-clique-relaxation});
$\SC = \{x\in [0,1]^{n+m} \mid \sum_{i\in \bar K_j}x_i=1,\ j\in[m]\}$ is a feasible set of the relaxed problem~\eqref{equ:clique-relaxation}~(\cref{sec:feasible-primal-estimate}).

\section{Overview of the known results and related work}
\label{sec:overview}

The MWIS problem is also known as \emph{maximum weight stable set} or \emph{vertex packing} and is mutually reducible to the \emph{minimum weight vertex cover}, \emph{maximum weight clique} and \emph{weighted set packing} problems~\cite{hochbaum1996approximating,eremeev2014optimal}.
Algorithms addressing these problems can be \emph{straightforwardly} applied to MWIS and vice versa.
Many theoretical results are shared by these problems as well.
Yet, we will stick to the MWIS problem formulation in the following.

\subsection{ILP formulation and edge relaxation}
Let $[n]$ denote $\{1, 2, \dots, n\}$ and assume that the set of graph nodes consists of $n$ elements, \ie, $\SV = [n]$.
We introduce a binary variable $x_i \in \{0, 1\}$ for each node $i \in \SV$.
If $x_i = 1$, it indicates that node $i$ is part of the maximum-weight independent set.
This results in a natural integer linear program (ILP) formulation of the MWIS problem:
\begin{equation}\label{equ:MWIS-edge}
   \max_{\vx\in\{0,1\}^n}\langle \vc,\vx\rangle
   \qquad \text{s.t.} \qquad
   x_i+x_j\le 1\,,\ \{i,j\}\in\SE \,.
\end{equation}
The respective LP relaxation, obtained by substituting the integer constraints $\vx\in\{0,1\}^n$ with the box constraints $\vx\in [0,1]^n$, is referred to as the \emph{edge relaxation}.
The edge relaxation is well-known to be tight for bipartite graphs~\cite{nemhauser1975vertex} and can be reduced to the minimum s-t cut problem (or to equivalent maximum flow), allowing it to be efficiently
solved~\cite{nemhauser1975vertex,hochbaum1996approximating}.
The efficient solvability of the edge relaxation follows from its two interrelated properties~\cite{nemhauser1975vertex}:
First, its corresponding polytope has only \emph{half-integral} vertices, meaning that all basic solutions of the relaxation have coordinates in the set $\{0, \nicefrac{1}{2}, 1\}$.
Second, it exhibits the \emph{persistency} property:
All coordinates of a relaxed solution with integer values (0 or 1) belong to some~(but the same) optimal solution of the non-relaxed problem~\eqref{equ:MWIS-edge}.
This persistency property has been particularly useful in improving approximation bounds for several primal heuristics, as demonstrated, \eg, by \cite{hochbaum1996approximating,kako2005approximation}.

However, the edge relaxation has its limitations: It gets loose for dense graphs, as illustrated by
\begin{proposition}\label{prop:edge-LP-solution-for-complete-graph}
  Let $G$ be fully-connected, all costs positive, $c_j > 0$, $j\in\SV$, and $c_i < \sum_{j\in\SV\backslash\{i\}}c_j$ for $i=\argmax_{j\in \SV}c_j$.
  Then the edge relaxation has a unique solution $(\underbrace{\nicefrac{1}{2},\nicefrac{1}{2},\dots,\nicefrac{1}{2}}_{n})$.
\end{proposition}
The statement follows from the half-integrality of the relaxation polytope, see \cref{sec:proof:edge-LP-solution-for-complete-graph} for the proof.

Note that cost positivity can be assumed \WLOG\ for the MWIS problem.
Nodes with non-positive costs can be excluded from the graph together with their incident edges without influencing the optimal objective value.

\Cref{prop:edge-LP-solution-for-complete-graph} shows that the edge relaxation may not reveal \emph{any} information about optimal solutions of the non-relaxed problem~\eqref{equ:MWIS-edge} for dense graphs.
To address this limitation, the \emph{clique constraints} and the respective \emph{clique relaxation} are considered for the MWIS problem in the literature~\cite{rodriguez2019persistency}, \cite[Sec.~64]{schrijver2003combinatorial}.
Before defining it, we first introduce a family of relaxations that includes both the edge and the clique relaxations.

\subsection{Clique cover relaxation family}\label{sec:clique-cover-family}
Let $K_j\subseteq\SV$, $j\in [m]$, be subsets of nodes such, that the respective induced subgraphs of $\SG$ are cliques.
Abusing notation, we will refer to $K_j$ themselves as \emph{cliques}.
We assume that the sets $K_j$, $j \in [m]$, cover all edges and nodes of the graph, \ie, $\SE = \bigcup_{j=1}^{m} \{\{i, i'\} \mid i, i' \in K_j\}$ and $\SV = \bigcup_{j=1}^{m} \{i \mid i \in K_j\}$, while explicit node coverage is required for isolated nodes. Further, we refer to $\{K_j\}_{j=1}^m$ as \emph{clique cover}.
We define the \emph{clique cover ILP formulation} of the MWIS problem as
\begin{equation}\label{equ:MWIS-clique}
   \max_{\vx\in\{0,1\}^n}\langle \vc,\vx\rangle
   \qquad \text{s.t.} \qquad
   \sum_{i\in K_j}x_i \le 1\,, j\in[m]\,.
\end{equation}
Its natural LP relaxation, further referred to as \emph{clique cover relaxation} is obtained by substituting integer constraints $\vx\in\{0,1\}^n$ with the respective box constraints $\vx\in [0,1]^n$.
The relaxation of~\eqref{equ:MWIS-clique} is called \emph{clique relaxation}, if the set $\{K_j \mid j\in [m]\}$, consists of \emph{all} cliques of $G$.

Whereas the solutions of the non-relaxed problems~\eqref{equ:MWIS-edge} and~\eqref{equ:MWIS-clique} coincide, their relaxations may differ significantly.
In particular, the clique relaxation lacks both the half-integrality and persistency properties in general~\cite{rodriguez2019persistency}.
Furthermore, it does not admit a reduction to any efficiently solvable LP subclass like the min-cost-flow problem. This is because its relaxation is as difficult as linear programs in general~\cite{prusa2019solving}.
However, it has the advantage of being a tighter relaxation in general, as shown by the following statement:
\begin{proposition}\label{prop:tightness-of-clique-relaxations}
  Let $K_j$, $j\in [m]$, and $K'_j$, $j\in [m']$, be two clique covers of the graph $G$ such that for every $K'_j$, $j\in [m']$, there is $l\in [m]\colon K'_j\subseteq K_l$.
  Then the relaxation determined by $K_j$, $j\in [m]$, is at least as tight as the relaxation determined by $K'_j$, $j\in [m']$.
\end{proposition}
\begin{proof}
  The proof is a direct consequence of the fact that for any $\vx\in[0,1]^n$ an inequality
  $\sum_{i\in K_l}x_i \le 1$ implies $\sum_{i\in K'_j}x_i \le 1$ as long as $K'_j\subseteq K_l$.
\end{proof}
\cref{prop:tightness-of-clique-relaxations} defines a partial order on the set of clique cover relaxations with respect to their tightness.
The edge relaxation is the least tight and constitutes the minimum with respect to this partial order.
The maximum is determined by the clique relaxation.
According to \cref{prop:tightness-of-clique-relaxations}, the latter is equivalent to the relaxation, where $K_j$, $j\in [m]$, consists of all maximal cliques (\ie, those that are not subgraphs of any other clique) only.
It is also known that clique inequalities are facet-defining for the stable set polytope, \ie, convex hull of all feasible solutions of~\eqref{equ:MWIS-edge},
if and only if the respective clique is maximal~\cite[Thm.~64.5]{schrijver2003combinatorial}.
In contrast to the edge relaxation, the maximal clique relaxation is also tight for complete graphs, see \cref{prop:clique-LP-solution-for-complete-graph} in \cref{sec:prop:clique-LP-solution-for-complete-graph}.
Unfortunately, the number of maximal cliques may grow exponentially with the graph size in general~\cite{moon1965cliques}, although they can be enumerated efficiently \cite{bron1973algorithm}.
Hence, in practice, one has to ideally consider a clique cover consisting of a subset of the maximal cliques. Such a cover can be built with a greedy algorithm that iteratively picks an uncovered node (or edge), grows it greedily to a maximal clique, adds this clique to the cover and marks covered nodes (edges).
Nevertheless, the algorithms we describe in this paper do \emph{not} require cliques in the cover to be maximal.

\subsection{Equality formulation}
Note that the clique constraints in~\eqref{equ:MWIS-clique} can be turned to equality by introducing $m$ \emph{slack variables} $x_{n+j}$, $j\in [m]$, and assuming $c_{n+j}=0$ and $\bar K_j=K_j\cup \{n+j\}$.
The respective LP relaxation takes the form:
\begin{equation}\label{equ:clique-relaxation}
   \max_{\vx\in [0,1]^{n+m}}\langle \vc,\vx\rangle
   \qquad \text{s.t.} \qquad
   \sum_{i\in \bar K_j}x_i = 1\,, j\in[m]\,.
\end{equation}
As~\eqref{equ:clique-relaxation} is equivalent to the clique cover relaxation, we will primarily use this name to refer to~\eqref{equ:clique-relaxation} unless otherwise stated.
Likewise, the non-trivial constraints of~\eqref{equ:clique-relaxation} will be called \emph{clique constraints} just like those in~\eqref{equ:MWIS-clique}.
To avoid ambiguity we will reference the specific formula.
\begin{remark}\label{rem:box-constraint-eq-ge-0}
  The box constraints $\vx\in[0,1]^{n+m}$ in~\eqref{equ:clique-relaxation} can be substituted with the non-negativity constraints $\vx\ge \vZero$ without changing the feasible set of the problem.
  This is because $\vx\le \vOne$ due to the clique constraints in~\eqref{equ:clique-relaxation}, as cliques cover all vertices of the underlying graph.
  In the following we stick to the box constraints for technical reasons.
\end{remark}

\subsection{Relation to pseudo-boolean optimization and discrete graphical models}\label{sec:relation-to-QPBO}
There is a natural relation between MWIS and pseudo-boolean optimization, see~\cite[Thm.~3, 4]{boros2002pseudo}. In particular, the MWIS problem~\eqref{equ:MWIS-edge} reduces
to the maximization of the \emph{quadratic} pseudo-boolean function
\begin{equation}
  \sum_{i=1}^{n} c_i \, x_i - M \cdot \sum_{\{i,j\}\in\SE}x_i x_j
\end{equation}
for sufficiently large $M > 0$.
The pseudo-boolean functions, in turn, have a one-to-one correspondence to the \emph{binary pairwise discrete energy minimization}~\cite[Ch.~12]{savchynskyy2019discrete}.
The respective \emph{local polytope} relaxation~\cite[Ch.~4]{savchynskyy2019discrete} is equivalent to the edge relaxation~\cite{sanghavi2007message}.
Some MWIS algorithms therefore adopt techniques from the field of discrete graphical models~\cite{shah2005max,sanghavi2007message}, however, their guaranteed performance does not go beyond the edge relaxation.

\subsection{Solvability and approximation}
The MWIS problem is polynomially solvable for perfect graphs~\cite{grotschel1984polynomial}.
In particular, this applies to bipartite graphs mentioned above.
Efficient methods also exist for other graph families, see~\cite[Ch.~64.9a]{schrijver2003combinatorial} for an overview.
Contrary to the closely related vertex cover problem, there are no constant factor approximation algorithms for MWIS~\cite{trevisan2014inapproximability}.
A number of approximation bounds has been obtained for the general problem as well as for its special cases~\cite{hochbaum1996approximating,paschos1997survey}, especially \wrt\ to the maximum degree of the underlying graph~\cite{kako2005approximation}.
These algorithms often make use of the persistency property of the edge relaxation.

\subsection{Primal heuristics}
Primal heuristics are a popular class of algorithms addressing MWIS, see \eg \cite{pullan2009optimisation,brendel2010segmentation,nogueira2018hybrid,dong2025metaheuristic,grossmann2025concurrent} and references therein.
Apart from the direct application, they can be used in combination with the relaxation-based methods, as in~\cite{dong2025metaheuristic} and this work, or, \eg, within reduction-based methods~\cite{lamm2019exactly}.
We get back to the first case in \cref{sec:primal-heuristics} and empirically compare to algorithms of both types in \cref{sec:experiments}.

\subsection{Relation to linear assignment}
The MWIS problem~\eqref{equ:MWIS-clique} can be seen as a generalization of the \emph{incomplete} linear assignment problem~\cite{dlask2025relative}, as the ILP formulation of the latter, namely
\begin{equation}\label{equ:LAP-def}
  \max_{\vx\in\{0,1\}^{n\times m}}\langle \vc,\vx\rangle
  \qquad \text{s.t.} \quad
  \sum_{i\in [n]}x_{ij} \le 1\,, j\in[m]
  \quad\text{and}\quad
  \sum_{j\in [m]}x_{ij} \le 1\,, i\in[n]\,,
\end{equation}
is a special case of~\eqref{equ:MWIS-clique}.
Incomplete linear assignment is tightly related and reducible to the (complete) linear assignment problem~\cite{haller2022comparative,dlask2025relative} in linear time.
From the other viewpoint, linear assignment is a specialization of the optimal transport problem.
These facts are important, as the Sinkhorn algorithm and its variants~\cite{peyre2017computational} have shown their remarkable efficiency and scalability for the optimal transport problems.
The Sinkhorn algorithm itself~\cite{sinkhorn1964relationship} is a special case of the famous Bregman method for linear programs~\cite{bregman1967relaxation}, which is one of the widely considered entropy optimization techniques~\cite{fang2012entropy}.
This leads us to the main contribution of our work formulated below.

\section{Contribution and paper content}

We introduce a scalable iterative algorithm for the clique cover family of LP relaxations~\eqref{equ:clique-relaxation}, based on the Bregman method -- interpretable as smoothed dual coordinate minimization and a generalization of the Sinkhorn algorithm for optimal transport. We propose several algorithmic variants that differ in numerics, smoothing schedule, and sparsity-based speedups.
As a dual method, our algorithm operates over the feasible set of the dual problem, while the corresponding primal variables remain infeasible until convergence. To address this, we propose an efficient projection scheme that maps dual iterates to the feasible set of the primal relaxed problem, enabling accurate duality gap estimation. This, in turn, allows for a novel duality-gap-based smoothing schedule, which consistently outperforms standard feasibility-based schemes in our experiments.

Our algorithm inherits convergence guarantees from the Bregman method and empirically achieves high-quality solutions. On large-scale instances (over 100k graph edges), it notably outperforms simplex and interior-point methods in runtime.

We further combine our dual solver with a simple yet effective primal heuristic to tackle the non-relaxed MWIS problem. Experiments show that our method produces high-quality solutions significantly faster than Gurobi~\cite{gurobi}, reduction-based approaches~\cite{lamm2019exactly}, and the recent meta-heuristic METAMIS~\cite{dong2025metaheuristic}, outperforming the latter in both speed and objective value.

Our code is publicly available on \href{https://github.com/vislearn/libmpopt/}{\texttt{github.com/vislearn/libmpopt/}}.

\subsection{Paper content}
We analyze the clique cover relaxation~\eqref{equ:clique-relaxation} and its Lagrange dual in~\cref{sec:dual-of-clique-relaxation}, starting with a natural baseline: non-smoothed coordinate minimization. While intuitive, it lacks convergence guarantees to the dual optimum.

In~\cref{sec:Bregman-for-MWIS}, we introduce a Bregman-based smoothed variant of coordinate minimization tailored to MWIS. The smoothing is controlled by a \emph{temperature} parameter: Higher temperatures yield faster convergence but coarser approximations.

\Cref{sec:numerics} presents two numerically stable implementations of our Bregman variant, adapted from similar techniques in optimal transport.

\Cref{sec:T-scheduling} discusses two strategies of smoothing scheduling and introduces a projection step that recovers feasible primal estimates, enabling duality gap evaluation.

To accelerate convergence, \cref{sec:sparsity-based-speedup} proposes two sparsity-based strategies: one heuristic, the other based on upper-bounding changes to the smoothed dual objective.

\Cref{sec:primal-heuristics} details a simple primal heuristic for constructing integer MWIS solutions. While basic, it performs well when combined with our dual solver.

Finally, \cref{sec:experiments} presents evaluations on two real-world benchmarks, related to image segmentation~\cite{prakash2022fully} and vehicle routing~\cite{dong2021new}. We study variants of our method and compare against Gurobi, reduction-based methods~\cite{lamm2019exactly,kamis-code}, and heuristics~\cite{nogueira2018hybrid,dong2025metaheuristic}.

\section{Dual of the clique cover relaxation}\label{sec:dual-of-clique-relaxation}

\paragraph{Dual problem.}
First, we introduce the constraint sets $J_i=\{j\in[m] \mid \bar K_j\ni i\}$ as the set of cliques containing the variable $x_i$.
In particular, for any slack variable $x_{n+j}$ it holds $J_{n+j}=\{j\}$.

By dualizing the equality constraints of the clique cover relaxation~\eqref{equ:clique-relaxation} we obtain its Lagrange dual problem:
\begin{equation}
  \label{equ:MWIS-dual-dummy}
   \min_{\vlambda\in\mathbb R^{m}}
    \Biggl[
      \;
      \max_{\vx\in[0,1]^{n+m}} \langle \vc,\vx\rangle
      - \sum_{j=1}^{m}\lambda_j\biggl(\sum_{i\in\bar K_j}x_i-1\biggr)
      \;
    \Biggr]
  =  \min_{\vlambda\in\mathbb R^{m}}
    \Biggl[
      \;
      \underbrace{
        \sum_{j=1}^{m}\lambda_j + \max_{\vx\in[0,1]^{n+m}}\langle \vc^\vlambda,\vx\rangle
      }_{D(\vlambda)}
      \;
    \Biggr]
\end{equation}
where $\vlambda$ are dual variables and $\vc^\vlambda_i=c_i - \sum_{j\in J_i}\lambda_j$ are the \emph{reduced} or \emph{reparametrized} costs.
From the first term in~\eqref{equ:MWIS-dual-dummy} it follows that $\langle\vc,\vx\rangle=\langle\vc^\vlambda,\vx\rangle + \sum_{j=1}^{m}\lambda_j$ holds for any $\vx$ that satisfies the clique constraints in~\eqref{equ:clique-relaxation}.
In other words, the costs of all feasible solutions (that includes integral solutions) of the problem~\eqref{equ:clique-relaxation}, are shifted by the constant $\sum_{j=1}^{m}\lambda_j$ with the transformation $\vc\to\vc^\vlambda$.
In particular, it implies that the set of optimal solutions of these problems is invariant with respect to this transformation.
Therefore, the latter is sometimes referred to as \emph{equivalent transformation} in the literature.

The dual objective $D(\vlambda)$ is convex piecewise affine as a maximum of a finite number of affine functions, since $\max_{\vx\in[0,1]^{n+m}}\langle \vc^\vlambda,\vx\rangle=\max_{\vx\in\{0,1\}^{n+m}}\langle \vc^\vlambda,\vx\rangle$ in~\eqref{equ:MWIS-dual-dummy}.
By strong duality, the minimum of the dual problem~\eqref{equ:MWIS-dual-dummy} coincides with the maximum of the clique relaxation~\eqref{equ:clique-relaxation}.

We introduce the set $\Sign(\vc) \coloneq \argmax_{\vx\in [0,1]^{n+m}}\langle \vc,\vx\rangle$.
Hence, $\vx\in\Sign(\vc)$ if
\begin{equation}\label{equ:MWIS-primal-rounded-solution}
  x_i =
  \begin{cases}
    0,         & \text{if}\ c_i < 0\,, \\
    1,         & \text{if}\ c_i > 0\,, \\
    \in [0,1], & \text{if}\ c_i = 0\,.
  \end{cases}
\end{equation}
The dual value is therefore equal to $D(\vlambda)=\sum_{j=1}^{m}\lambda_j + \langle \vc^\vlambda,\vx^*\rangle$ for any $\vx^*\in\Sign(\vc^{\vlambda})$.
Different subgradients of $D$ can be obtained for different values of $\vx^*\in\Sign(\vc^{\vlambda})$ as
\begin{equation}\label{equ:MWIS-dual-subgradient}
  \frac{\partial D}{\partial \lambda_j}[\vx^*]
  = 1 + \frac{\partial \langle \vc^\vlambda,\vx^*\rangle}{\partial \lambda_j}
  = 1 +  \Bigl\langle \frac{\partial \vc^\vlambda}{\partial \lambda_j},\vx^*\Bigr\rangle
  = 1 - \sum_{i\in \bar K_j}x_i^*\,.
\end{equation}
As follows from~\eqref{equ:MWIS-dual-subgradient}, the considered subgradient coordinate $\frac{\partial D}{\partial \lambda_j}[\vx^*]$ is equal to $0$ if the sum of non-zero coordinates of $\vx^*$ in the clique $\bar K_j$ is equal to $1$.
This implies that $\vlambda$ is dual optimal if and only if $\Sign(\vc^\vlambda)$ contains at least one vector feasible for~\eqref{equ:clique-relaxation}.
This vector is a solution of the relaxed primal problem:
\begin{proposition}\label{prop:dual-optimality-condition}
  Let $\vx^*\in\Sign(\vc^\vlambda)$ and $\frac{\partial D}{\partial \lambda_j}[\vx^*]=0$ for all $j\in[m]$.
  Then $\vx^*$ is a solution of the respective clique relaxation~\eqref{equ:clique-relaxation}.
  If additionally $\vx^*\in\{0,1\}^{n+m}$, then its first $n$ coordinates constitute a solution of the non-relaxed problem~\eqref{equ:MWIS-clique}.
\end{proposition}
\begin{proof}
  The condition $\frac{\partial D}{\partial \lambda_j}[\vx^*]=0$ for all $j\in[m]$ together with~\eqref{equ:MWIS-dual-subgradient} implies feasibility of $\vx^*$ and optimality of $\vlambda$. Therefore, for any $\vx$ feasible for~\eqref{equ:clique-relaxation} it holds
  $
    \langle \vc,\vx\rangle \le D(\vlambda) = \sum_{j=1}^{m}\lambda_j + \langle \vc^\vlambda,\vx^*\rangle = \langle \vc,\vx^*\rangle
  $,
  and thus, $\vx^*$ is an optimal solution of~\eqref{equ:clique-relaxation}.
  If additionally $\vx^*$ is integer, the inequality proves its optimality for the non-relaxed problem~\eqref{equ:MWIS-clique} as well.
\end{proof}

\paragraphRunIn{The dual coordinate minimization \cref{alg:non-smooth-bcd}} iteratively considers each clique and changes the value of the respective dual variable.
This is done such that the largest reduced cost of the nodes belonging to the clique becomes zero.

\begin{algorithm}[t]
  \addtolength{\hsize}{1.5em}%
  \DontPrintSemicolon
  \KwIn{$\vlambda$ \tcc*[h]{initial dual solution}} 
  \While{stopping criterion not fulfilled}{
    \For(\tcc*[h]{loop over cliques}){$j\in [m]$}{
      \tcc{turn the max.~cost in the clique to $0$:}
      $i^*\in\argmax_{i\in \bar K_j}c^{\vlambda}_i$ \tcc{pick the clique-largest reduced cost}
      $\lambda_j \coloneq \lambda_j+c^{\vlambda}_{i^*}$ \tcc{\dots\ and make it zero}
    }
  }
  \KwResult{$\vlambda$}
  \caption{Dual coordinate descent algorithm.}
  \label{alg:non-smooth-bcd}
\end{algorithm}

Being very simple, this algorithm implements the \emph{coordinate minimization} principle.
Indeed, the dual objective $D$ is convex.
Therefore, its restriction to the $j$-th variable $\lambda_j$ is convex as well.
In turn, a necessary and sufficient optimality condition for the restricted convex function is existence of a zero subgradient.
In this case, after each inner iteration of \cref{alg:non-smooth-bcd} there should exist $\vx^*\in\Sign(\vc^\vlambda)$ such that $\frac{\partial D}{\partial \lambda_j}=1 - \sum_{i\in \bar K_j}x_i^*=0$ (see~\eqref{equ:MWIS-dual-subgradient}) for the value of $j$ being an index of the currently processed clique.
To obtain such $\vx^*$ it is sufficient to set $x^*_{i^*}=1$ and $x^*_i=0$ for $i\in \bar K_j\backslash\{i^*\}$, where $i^*$ is defined as in \cref{alg:non-smooth-bcd}.

Unfortunately, coordinate minimization is not guaranteed to converge to the optimum for convex piece-wise affine functions~\cite[Ch.~2.7]{bertsekas1999nonlinear}.
This holds in general, as well as for the dual~\eqref{equ:MWIS-dual-dummy} in particular,
see \cref{ex:bcd-fix-point} in \cref{sec:ex:bcd-fix-point}.
Consequently, the fix-points attained by \cref{alg:non-smooth-bcd} in our experiments give only a loose approximation of the dual optimum, although the respective lower bounds are usually much better than those of the edge relaxation, see \cref{fig:plots-relaxed}(a) in \cref{sec:experiments}.

\section{Bregman method for MWIS}\label{sec:Bregman-for-MWIS}

The method of~\cite{bregman1967relaxation} can be used to address linear programs $\{\max_{\vx\ge 0}\langle \vc,\vx\rangle\,,\text{ s.t. } A\vx=\vb\}$ through solving their approximations
$\{\max_{\vx\ge 0}\langle \vc,\vx\rangle + T \, \SH(\vx)\,,\ \text{s.t.}\ A\vx=\vb\}$ with the function $\SH(\vx) \coloneq - \sum_{i=1}^{n+m}\left(x_i\log x_i-x_i\right)$ being the entropy shifted by a linear function for technical reasons. We refer to~\cite{bregman1967relaxation} for a more general treatment.
Parameter $T > 0$ is the \emph{temperature} or \emph{smoothing} parameter that regulates the trade-off between approximation accuracy and convergence speed.
The problem above takes the form of the \emph{smoothed clique cover relaxation} in our case:
\begin{equation}\label{equ:smoothed-clique-relaxation}
   \max_{\vx\in [0,1]^{n+m}}\langle \vc,\vx\rangle  + T \, \SH(\vx)
   \qquad \text{s.t.} \qquad
   \sum_{i\in \bar K_j}x_i = 1\,, j\in[m]\,,
\end{equation}
since the box constraints $\vx\in [0,1]^{n+m}$ in~\eqref{equ:smoothed-clique-relaxation} are equivalent to $\vx\ge \vZero$ as noted in \cref{rem:box-constraint-eq-ge-0}.

On each iteration of the Bregman method the current iterate $\vx$ is projected to the hyperplane determined by a single row of the constraint matrix $A$.
In case of the function $\SH$ defined as above and constrains of the form $\sum_{i\in \bar K_j}x_i=1$ as in~\eqref{equ:smoothed-clique-relaxation}, this projection results in \emph{normalization} of the respective coordinates of $\vx$, as given by \cref{alg:naive-sinkhorn}.
The only computation that depends on the cost vector $\vc$ is the initialization of the algorithm: Whereas it converges to a feasible point for \emph{arbitrary} initialization $\vx > \vZero$, its limit point is a solution to the problem~\eqref{equ:smoothed-clique-relaxation} only if properly initialized.
When applied to the linear assignment problem \cref{alg:naive-sinkhorn} is usually referred to as \emph{Sinkhorn algorithm} due to the seminal paper by \cite{sinkhorn1964relationship} showing that iterative normalization of a square matrix turns it into a bi-stochastic one.

\begin{algorithm}[t]
  \addtolength{\hsize}{1.5em}%
  \DontPrintSemicolon
  \KwIn{$\vlambda$, $T$ \tcc*[h]{initial dual solution, temperature}}
  \textbf{Init:} $x_i=\exp(\nicefrac{c^{\vlambda}_i}{T}),\ i\in [n]$\;
  \While{stopping criterion not fulfilled}{
    \For(\tcc*[h]{loop over cliques}){$j\in [m]$}{
      $s \coloneq \sum_{i\in \bar K_j}x_i$
      \tcc{compute the sum of $x_i$ within the clique}
      $x_i \coloneq x_i/s,\ \forall i\in \bar K_j$
      \tcc{\dots\ and make it equal to $1$} \nllabel{step:alg:naive-sinkhorn:update}
    }
  }
  \KwResult{$\vx$}
  \caption{Naive Bregman method for MWIS.}
  \label{alg:naive-sinkhorn}
\end{algorithm}

The following statements formalize the main properties of \cref{alg:naive-sinkhorn}:
\begin{lemma}\label{prop:x-le-0}
  After the first outer iteration of \cref{alg:naive-sinkhorn} it holds $\vZero \le \vx\le \vOne$ .
\end{lemma}
\begin{proof}
  The proof directly follows from the update in~\cref{step:alg:naive-sinkhorn:update} with $s\ge x_i$.
\end{proof}

\begin{proposition}\label{prop:naive-sinkhorn-convergence}
  Iterates of \cref{alg:naive-sinkhorn} converge to the solution of the {smoo\-thed} relaxation~\eqref{equ:smoothed-clique-relaxation}.
\end{proposition}
\begin{proof}
\newcommand\MyBregman[1]{\cite[#1]{bregman1967relaxation}}
Denoting $S\bar \coloneq \BR^{n+m}_{\ge 0}$ and taking into account \cref{rem:box-constraint-eq-ge-0}, we conclude that~\eqref{equ:smoothed-clique-relaxation} is equivalent to
  $
    \left\{
      \min_{\vx\in\bar S} \, \langle -\nicefrac{\vc}{T}, \vx \rangle - \SH(\vx) \mid A\vx=\vb
    \right\}
  $,
which is the problem format considered in \MyBregman{(2.1)--(2.3)}.
Here the matrix $A$ is defined by its rows $\va_{j}\in\BR^{n+m}$, $a_{ji}=\llbracket i\in \bar K_j\rrbracket$ and $b_j=1$ for $j\in[m]$.
The sufficient convergence conditions for the objective $f(\vx) \coloneq \langle -\nicefrac{\vc}{T}, \vx \rangle - \SH(\vx)$ are determined by its Bregman divergence $D_f(\vx,\vy)$ defined as in \MyBregman{(1.4)}.
As shown in \MyBregman{\textsection 1}, these conditions are satisfied for $D_g$ with $g(\vx)=\sum_{i=1}^{n+m}x_i\log x_i$, and $D_f(\vx,\vy)=D_g(\vx,\vy)$, as $f(\vx)-g(\vx)$ is an affine function. The latter implies that the conditions hold also for $f(\vx)$.

As stated in \MyBregman{\textsection 2}, an initial point $\vx^0$ of the algorithm must satisfy
$f'(\vx^0)=\vu A$ for some $\vu\in\BR^{m}$.
This translates into $\left(\langle -\nicefrac{\vc}{T}, \vx^0 \rangle - \SH(\vx^0)\right)' = -\nicefrac{\vc}{T} + \log \vx^0 = \vu A$. Solving this w.r.t. $\vx^0$, substituting $\vu \coloneq -\nicefrac{\vlambda}{T}$ and taking into account the form of $A$, we conclude that $\vx^0=\exp(\nicefrac{\vc^{\vlambda}}{T})$, which is the starting point of \cref{alg:naive-sinkhorn}.

According to \MyBregman{(2.7)--(2.8)}, the update rule for the $j$-th constraint can be derived from the equations $f'(\vx')=f'(\vx)+\gamma \va_j$ and $\sum_{i=1}^{n+m} a_{ji}x'_i=b_j$ for some $\gamma\in\BR$.
The first equation takes the form
$
  x'_i \coloneq x_i \exp(\nicefrac{\gamma a_{ji}}{T})
$, with $\gamma$ being computed as a solution of the second equation taking the form
$
  \sum_{i=1}^{n+m}a_{ji}x_i \exp(\nicefrac{\gamma a_{ji}}{T})=b_j
$.
Plugging the values of $a_{ji}$ and $b_j$ into the above equations and solving the latter w.r.t.\ $\gamma$ gives the update $x'_i = \frac{x_i}{\sum_{i\in \bar K_j}x_i}$ as in \cref{step:alg:naive-sinkhorn:update} of \cref{alg:naive-sinkhorn}.
\end{proof}

To be able to apply \cref{alg:naive-sinkhorn} to the clique relaxation~\eqref{equ:clique-relaxation}, one has to define \emph{how to select the temperature $T$} as well as \emph{how to deal with situations when $s=\sum_{i\in \bar K_j}x_i$ is numerically equal to zero}.
The latter may happen due to computation of $x_i=\exp(\nicefrac{c^{\vlambda}_i}{T})$ for low temperatures and negative costs, and triggers an error during division to~$s$.
We cover the second question first in \cref{sec:numerics} and get back to the first one in \cref{sec:T-scheduling}.

\section{Numerics}\label{sec:numerics}

\paragraph{Log-domain numerics.}
Consider \cref{alg:log-sinkhorn}.
Note that it differs from \cref{alg:non-smooth-bcd} by essentially the single term $T\log\sum_{i\in \bar K_j}\exp\left((c^{\vlambda}_i-c^{\vlambda}_{i^*})/T\right)$ in the computation pipeline.
\begin{algorithm}[t]
  \addtolength{\hsize}{1.5em}%
  \DontPrintSemicolon
  \KwIn{$\vlambda$, $T$ \tcc*[h]{initial dual solution, temperature}}
  \While{stopping criterion not fulfilled}{
    \For(\tcc*[h]{loop over cliques}){$j\in [m]$}{
      \tcc{make the largest cost in the clique $\approx 0$}
      $i^*\in\argmax_{i\in \bar K_j}c^{\vlambda}_i$
      \tcc{pick the clique-largest reduced cost}
      $\lambda_j \coloneq \lambda_j+c^{\vlambda}_{i^*}+T\log\sum_{i\in \bar K_j}\exp{\frac{c^{\vlambda}_i-c^{\vlambda}_{i^*}}{T}}$ \nllabel{alg:log-sinkhorn:lambda-update}
      \tcc{\dots\ and make it $\approx 0$}
    }
  }
  \KwResult{$\vlambda$}
  \caption{Bregman method in log domain.}
  \label{alg:log-sinkhorn}
\end{algorithm}

The following presentation recapitulates the well-known result for this type of algorithms in the optimal transport domain,
see also \cref{sec:prop:log-sinkhorn-equiv-to-naive-sinkhorn} for the proof:
\begin{lemma}\label{prop:log-sinkhorn-equiv-to-naive-sinkhorn}
  \cref{alg:log-sinkhorn} is equivalent to \cref{alg:naive-sinkhorn} through the transformation $\vx=\exp(\nicefrac{\vc^{\vlambda}}{T})$.
\end{lemma}

Contrary to \cref{alg:naive-sinkhorn}, \cref{alg:log-sinkhorn} is numerically stable.
Its stability follows from $c^{\vlambda}_i-c^{\vlambda}_{i^*} \le 0$ and, therefore, $\exp\left((c^{\vlambda}_i-c^{\vlambda}_{i^*})/T\right)\in (0,1]$, and $c^{\vlambda}_{i^*}-c^{\vlambda}_{i^*} = 0$ implying $\sum_{i\in \bar K_j}\exp\left((c^{\vlambda}_i-c^{\vlambda}_{i^*})/T\right)>1$.
The latter is sufficient to have a numerically stable logarithm computation.
\begin{lemma}\label{lem:costs-are-negative}
  After the first iteration of \cref{alg:log-sinkhorn} the reduced costs remain non-positive, \ie, $\vc^\vlambda \le \vZero$.
\end{lemma}
  The proof follows directly from \cref{prop:x-le-0} and \cref{prop:log-sinkhorn-equiv-to-naive-sinkhorn}.
\begin{proposition}\label{prop:bregman-smoothed-bca-optimum}
  \cref{alg:log-sinkhorn} converges to the optimum of the \emph{dual smoothed clique cover relaxation} problem
  \begin{equation}\label{equ:smoothed-dual}
    \min_{\vlambda}\Biggl[D^T(\vlambda) \coloneq  \sum_{j=1}^{m}\lambda_j
      + \max_{\vx\in [0,1]^{n+m}}\Bigl(\langle \vc^\vlambda,\vx\rangle+T \, \SH(\vx)\Bigr)\Biggr],
  \end{equation}
  and implements the coordinate minimization principle.
  In particular, $D^T$ is convex differentiable and right after processing the clique $j\in[m]$ in \cref{alg:log-sinkhorn:lambda-update} it holds $\frac{\partial D^T}{\partial \lambda_j}=0$.
\end{proposition}
The proof of proposition is based on the following technical
\begin{lemma}\label{lem:smooth-dual-explicit-formula}
  For any $T\ge 0$ and $\vc \le 0$ it holds:
 $
    \max\limits_{\vx\in [0,1]^{n+m}}\Bigl(\langle \vc,\vx\rangle+T\,\SH(\vx)\Bigr)= T\sum\limits_{i=1}^{n+m} \exp(\nicefrac{c_i}{T})
$.
\end{lemma}
\begin{proof}
  Consider
  \begin{equation}\label{equ:Delta-def-0}
    \max_{\vx\in [0,1]^{n+m}}\langle \vc,\vx\rangle - T\sum_{i=1}^{n+m}\left(x_i\log x_i-x_i\right)
    = \sum_{i=1}^{n+m} \max_{x_i\in [0,1]} c_i x_i - T x_i\log x_i +Tx_i\,.
  \end{equation}

  The function $f(x) \coloneq cx - T x\log x +Tx$ of a scalar $x$ dependent on scalars $c\le 0$ and $T\ge 0$ is concave and therefore its optimum is defined by the following condition
  \begin{equation}
    0=\nabla f(x) = c-T \, \biggl(\log x + \frac{x}{x} -1 \biggr) = c - T \, \log x\,.
  \end{equation}
  It implies
$\log x = \nicefrac{c}{T}$ and $x = \exp(\nicefrac{c}{T})$.
  Since $c\le 0$, it holds $0\le x \le 1$, \ie it satisfies the box constraints on $x_i$ in~\eqref{equ:Delta-def-0}.
  The value of $f(x)$ in its maximum is computed as
  \begin{equation}
    \max_{x\in[0,1]} f(x)
    = c\cdot \exp(\nicefrac{c}{T}) - T\cdot \exp(\nicefrac{c}{T}) \cdot \nicefrac{c}{T} + T \exp(\nicefrac{c}{T})
    = T\cdot \exp(\nicefrac{c}{T})\,.
  \end{equation}
  Putting all together proves the statement of the lemma.
\end{proof}

\begin{proof}[\cref{prop:bregman-smoothed-bca-optimum}]
  The first part of the statement directly follows from the strong duality, \cref{prop:naive-sinkhorn-convergence,prop:log-sinkhorn-equiv-to-naive-sinkhorn}.
  The second part follows from~\cref{lem:smooth-dual-explicit-formula} that implies that after processing the clique $\bar K_j$ it holds
  \def\tmpStackRel{\text{(Lem.~\ref{prop:log-sinkhorn-equiv-to-naive-sinkhorn})}}
  \def\tmpEqWithout{\stackrel{\phantom{\tmpStackRel}}{=}}
  \def\tmpEqWith   {\stackrel{         \tmpStackRel }{=}}
    \begin{equation}
    \frac{\partial D^T}{\partial \lambda_j}
    = 1 + T\sum_{i=1}^{n+m}\frac{\partial (\exp(\nicefrac{c^{\vlambda}_i}{T}))}{\partial \lambda_j}
    = 1 -  \sum_{i\in\bar K_j}\exp(\nicefrac{c^{\vlambda}_i}{T})
    \tmpEqWith 1 -  \sum_{i\in\bar K_j}x_i
    \stackrel{\text{(Alg.~\ref{alg:naive-sinkhorn})}}{=} 0
  \end{equation}
\end{proof}

\Cref{prop:naive-sinkhorn-convergence,prop:bregman-smoothed-bca-optimum} imply that the dual iterates $\vlambda$ of \cref{alg:log-sinkhorn} determine the primal ones $\vx=\exp(\nicefrac{\vc^{\vlambda}}{T})$ and the latter converge to the solution of the smoothed clique relaxation~\eqref{equ:smoothed-clique-relaxation}.

Along with the advantages such as numerical stability and possibility to estimate an upper bound due to the explicit computations of the dual iterates, \cref{alg:log-sinkhorn} has a major disadvantage: It has high computational cost due to extensive exponentiations.
An alternative way to perform computations that allows to keep advantages of \cref{alg:log-sinkhorn} and avoid its disadvantages is reviewed in the following paragraph.

\paragraphRunIn{Exp-domain numerics with stabilization} is represented by \cref{alg:exp-domain-stabilized}.
It keeps track of the exponentiated version $\alpha_j \coloneq \exp(-\nicefrac{\lambda_j}{T})$, $j\in [m]$, of the dual variables.
Note that for any $j\in [m]$ and any $i\in\bar K_j$, it holds
\begin{equation}
  x_i = \exp\Biggl[\frac{c^\vlambda_i}{T}\Biggr]
  = \exp\Biggl[\frac{1}{T} \, \biggl(c_i-\sum_{\mathclap{j'\in J_i}}\lambda_{j'}\biggr)\Biggr]
  = \exp\Biggl[\frac{1}{T} \, \biggl(c_i-\sum_{\mathclap{j'\in J_i\backslash\{j\}}}\lambda_{j'}\biggr)\Biggr]
  \underbrace{\exp\Biggl[-\frac{\lambda_j}{T}\Biggr]}_{\alpha_j}\,.
\end{equation}
So, $x_i$, where $i\in\bar K_j$, are proportional to $\alpha_j$.
Hence, instead of normalizing $x_i$, we can update $\alpha_j$ by dividing it by the normalization factor $s$, see \cref{alg4:line4}.
The main advantage of this approach is accessibility of the current  dual solution as $(\vlambda-T\log\valpha)$.
Since computation of the normalizing factor~$s$ requires values of~$\vx$, we update $\vx$ directly as well, see \cref{alg4:line6}.

\begin{algorithm}[t]
  \addtolength{\hsize}{1.5em}%
  \DontPrintSemicolon
  \KwIn{$\vlambda$, $T$ \tcc*[h]{initial dual solution, temperature}}
  \textbf{Init:} $\alpha_j=1$, $j\in [m]$; $x_i=\exp(\nicefrac{c^{\vlambda}_i}{T}),\ i\in [n]$\;
  \While{stopping criterion not fulfilled}{
    \For(\tcc*[h]{for all cliques}){$j\in[m]$}{
      $s \coloneq \sum_{i\in \bar K_j}x_i$ \nllabel{alg4:line3}
      \tcc{compute sum over the clique}
      $\alpha_j \coloneq \alpha_j/s$ \nllabel{alg4:line4}
      \tcc{update exponentiated dual}
      \If{$\alpha_j + 1/\alpha_j < \tauStab$ \nllabel{alg4:line5}}{
        $x_i \coloneq x_i/s\,,\ \forall i\in \bar K_j$ \nllabel{alg4:line6}
        \tcc{if $\alpha_j \approx 1$ run naive Bregman method}
      }
      \Else{

        \For(\tcc*[h]{\textbf{stabilize:} for all cliques}){$j\in[m]$}{ \nllabel{alg4:line8}
          $\lambda_j \coloneq \lambda_j - T\log{\alpha_j}$
          \tcc{push $\valpha$ into $\vlambda$ to obtain a dual solution}
          $\alpha_j \coloneq 1$
          \tcc{reinitialize $\valpha$}
       }
       $x_i \coloneq \exp(\nicefrac{c^{\vlambda}_i}{T}),\ \forall i\in \bar K_j$ \nllabel{alg4:line10}
       \tcc{update $x$ in a stable way}
      }
    }
  }
  \For(\tcc*[h]{for all cliques}){$j\in[m]$}{
    $\lambda_j \coloneq \lambda_j - T\log{\alpha_j}$
    \tcc{push $\valpha$ into $\vlambda$ to obtain a dual solution}
  }
  \KwResult{$(\vx,\vlambda)$}
  \caption{Stabilized Bregman method in exp-domain.}
  \label{alg:exp-domain-stabilized}
\end{algorithm}

As long as $\valpha$ remains close enough to $\vOne$ (\cref{alg4:line5}, see \cref{sec:experiments} for the used values of the \emph{stabilization threshold} $\tauStab$), this implies stable computations for $\valpha$ and $\vx$ (\cref{alg4:line4,alg4:line6}).
Should $\valpha$ become very large or very small (\cref{alg4:line5}), the computations may get imprecise and, therefore, one switches to the expensive computations in the log-domain to stabilize them (\crefrange{alg4:line8}{alg4:line10}).

Note, however, that close to the optimum of the smoothed problem~\eqref{equ:smoothed-clique-relaxation}, and, therefore, feasibility of~$\vx$, computations in \cref{alg:exp-domain-stabilized} \emph{do not require stabilization} as $s\approx 1$ and, hence, $\valpha\approx \vOne$ in this case.
This emphasizes importance of the proper scheduling of the temperature parameter $T$ as the latter influences the distance to the smoothed optimum.
We consider this question in \cref{sec:T-scheduling}.

\section{Smoothing\,/\,Temperature scheduling}\label{sec:T-scheduling}

The lower the temperature~$T$ is, the better the smoothed clique cover relaxation~\eqref{equ:smoothed-clique-relaxation} approximates its non-smoothed counterpart~\eqref{equ:clique-relaxation}, see~\cite{weed2018explicit}.
Given a required approximation level, one could imagine two strategies of temperature selection.
The first one is the \emph{fixed} temperature.
Although being very simple, it leads to very slow algorithms.
First, because the Bregman algorithm makes only small progress for low temperatures.
Second, because the numerical stabilization in \cref{alg:exp-domain-stabilized} can be required very often for computations far from the optimum.
This would eliminate the speed advantage of the exp-domain-based \cref{alg:exp-domain-stabilized} over the log-domain-based \cref{alg:log-sinkhorn}.

Hence, a practical solution is the \emph{temperature scheduling}, \ie, starting with a high value of $T$ and decreasing it as the computation progresses.
We consider two methods for the temperature scheduling below.

\subsection{Method~1: Feasibility-based}\label{sec:feasibility-constraint-satisfaction}
Feasibility constraints satisfaction is the standard stopping criterion in the optimal transport literature~\cite{schmitzer2019stabilized}, used to decide whether it is time to decrease the temperature.
Typically, one checks whether
\begin{equation}\label{equ:eps-constraint-satisfaction}
  \Biggl|\sum_{i\in \bar K_j}x_i - 1\Biggr| \le \tauFeas\,, \forall j\in[m]\,,
\end{equation}
for some \emph{feasibility threshold} $\tauFeas\ge 0$.
Should the constraints~\eqref{equ:eps-constraint-satisfaction} be fulfilled, one drops the temperature by the predefined factor $0<\tauDrop<1$, \ie ${T \coloneq \tauDrop \cdot T}$.
For higher efficiency the conditions~\eqref{equ:eps-constraint-satisfaction} are checked only each $\tauBatch\in\BN$ iterations.

\begin{remark}
  To simplify the selection of the initial temperature and keep computations numerically stable, we:
  \emph{(i)}~Run one iteration of \cref{alg:non-smooth-bcd} at the beginning to ensure $\vc^\vlambda < 0$, and
  \emph{(ii)}~scale all costs such that they lie in the interval $[-1,0]$.
  After completing these operations, we launch \cref{alg:exp-domain-stabilized} with an initial temperature of $T=0.01$.
\end{remark}

\begin{remark}
  There is an alternative method for feasibility-based temperature scheduling, introduced by~\cite{kushinsky2019sinkhorn} for the quadratic assignment problem and equivalent to the above rule with $\tauDrop=0.5$.
  Although authors report its numerical stability in their experiments, for other problem instances it may require stabilization as well.
\end{remark}

\subsection{Method~2:  Duality gap-based}\label{sec:duality-gap-scheduling}
Duality gap is an alternative way to control the distance to the optimum when dual variables $\vlambda$ are computed explicitly, \ie, in \cref{alg:log-sinkhorn,alg:exp-domain-stabilized}.
Let $\hat \vx$ and $\vx^*$ be a \emph{feasible} and an \emph{optimal} solutions to the smoothed clique cover relaxation~\eqref{equ:smoothed-clique-relaxation} respectively.
Let also $\vlambda$ and $\vlambda^*$ be a \emph{feasible} and an \emph{optimal} solutions to its smoothed dual~\eqref{equ:smoothed-dual}.
The strong duality implies $\langle \vc,\hat \vx\rangle+T \, \SH(\hat \vx)\le \langle \vc,\vx^*\rangle+T \, \SH(\vx^*)=D^T(\vlambda^*) \le D^T(\vlambda)$.
A progress of the Bregman method is assured if the inequalities above hold strictly, \ie, $T < \frac{D^T(\vlambda)-\langle \vc,\hat \vx\rangle}{\SH(\hat \vx)}$.
This inequality turns into equality in the optimum of the smoothed problems~\eqref{equ:smoothed-clique-relaxation}-\eqref{equ:smoothed-dual}.
Hence, the respective temperature update condition, for a suitable predefined ${0<\tauGap < 1}$, may look like
\begin{equation}\label{equ:duality-gap-T-condition}
  \text{if}\quad T > \tauGap\frac{D^T(\vlambda)-\langle \vc,\hat \vx\rangle}{\SH(\hat \vx)} \quad \text{update}\ T\,.
\end{equation}

To ensure monotonicity of the temperature, we update it by setting to
\begin{equation}\label{equ:duality-gap-T-scheduling}
  T \coloneq \min\left\{T, \; \tauGap\frac{D^T(\vlambda)-\langle \vc,\hat \vx\rangle}{\SH(\hat \vx)}\right\} \end{equation} each $\tauBatch\in\BN$ iterations.
Between temperature changes, we use \cref{alg:log-sinkhorn}, or~\ref{alg:exp-domain-stabilized}, initialized with the current dual solution $\vlambda$.

\subsection{Feasible primal estimates for Method 2}\label{sec:feasible-primal-estimate}

\paragraph{Obtaining feasible primal estimates.}
Let $\SC$ stand for the feasible set of~\eqref{equ:smoothed-clique-relaxation}.
Computations~\eqref{equ:duality-gap-T-condition} and~\eqref{equ:duality-gap-T-scheduling} require a \emph{feasible} $\hat \vx$, \ie, $\hat \vx \in \SC$.
Contrary to the unconstrained dual iterates $\vlambda$, their primal counterparts $\vx=\exp(\nicefrac{\vc^{\vlambda}}{T})$ are infeasible until convergence, \ie, $\vx\notin \SC$.
To obtain feasible iterates we \emph{project} $\vx$ to the set $\SC$:
\begin{definition}
  A mapping $\SP\colon \BR^{n+m}_{\ge 0}\to\SC$ is called \emph{projection}, if it is \emph{idempotent}, \ie, $\SP(\SP(\vx))=\SP(\vx)$, and \emph{continuous}.
  The latter guarantees that $\vx^t\stackrel{t\to\infty}{\to} \vy\in \SC$ implies $\SP(\vx^t)\stackrel{t\to\infty}{\to} \SP(\vy)=\vy$.
\end{definition}
Consider \cref{alg:truncation-projection} that introduces what we call a \emph{truncation projection}.
The algorithm checks all coordinates of $\vx$ and fixes their values in the projection $\hat\vx$ as long as they do not violate the clique constraints.
The value is truncated if necessary to avoid constraint violation.
The slack variables $\hat x_{n+j}$ are set to guarantee equality in the clique constraints.

\begin{algorithm}[t]
  \addtolength{\hsize}{1.5em}%
  \DontPrintSemicolon
  \KwIn{$\vx\in \BR^{n+m}_{+}$}
  \tcc{initialize the feasible solution:}
  $\hat\vx\coloneq\mathbf{0}\in \BR^{n+m}$ \tcc{set all coordinates to $0$ and}
  $\hat x_{n+j}=1,\ j\in [m]$ \nllabel{alg5:line2} \tcc{set slack variables to $1$ for all $m$ cliques}

  \For(\tcc*[h]{update all non-slack variables}){$i\in [n]$ \nllabel{alg5:line3}}{
    $M \coloneq \min_{j\in J_i} \hat x_{n+j}$
    \tcc{max.~value of $x_i$ fitting clique constraints}
    $\hat x_i \coloneq \min\{x_i,M\}$ \nllabel{alg5:line5}
    \tcc{truncate $x_i$ if needed}
    \For{$j\in J_i$}{
      $\hat x_{n+j} \coloneq \hat x_{n+j} -  \hat x_i$ \nllabel{alg5:line7} \tcc{decrease the slacks resp.}
    }
  }
  \KwResult{$\SP(\vx) \coloneq \hat \vx$}
  \caption{Truncation projection.}
  \label{alg:truncation-projection}
\end{algorithm}

\Cref{alg:truncation-projection} has the following properties:
\begin{lemma}[Feasibility]\label{lem:truncation-projection-is-feasible}
  The vector $\hat\vx$ is always feasible at the end of each iteration of the outer loop defined in \cref{alg5:line3}.
\end{lemma}
\begin{proof}
  Initialization: $\hat \vx$ is feasible by construction in \cref{alg5:line2}.
  Exactly two coordinates in each inequality $j\in J_i$ containing $x_i$ are changing at each iteration of the loop over $i\in[n]$ in \crefrange{alg5:line3}{alg5:line7}: The value of $\hat x_i$ is increased (\cref{alg5:line5}) and the value of the respective slack $\hat x_{n+j}$ is decreased by the same amount.
  So the sum of coordinates of $\hat\vx$ in each clique constraint remains unchanged and equal to one.
\end{proof}
\begin{lemma}[Idempotence]\label{lem:truncation-projection-is-idempotent}
  If the input $\vx$ is feasible, then $\SP(\vx)=\vx$.
\end{lemma}
\begin{proof}
  As $\sum_{i'\in\bar K_j}x_{i'}=1$ for all $j\in J_i$, it holds $M\ge x_i$ in \cref{alg5:line5}.
  In turn it implies $\hat x_i=x_i$.
\end{proof}

\begin{lemma}[Continuity] \label{lem:truncation-projection-is-continuous}
  The mapping $\SP\colon \BR^{n+m}_{\ge 0}\to \SC$ is continuous.
  Moreover,
  \begin{equation}
    \biggl\vert \sum_{i\in \bar K_j}x_i - 1 \biggr\vert \le \epsilon\,,\ \forall j\in[m]
  \end{equation}
  implies
  \begin{equation}
    \Bigl\vert x_i - \SP(\vx)_i\Bigr\vert \le \epsilon,\ \forall i\in[n]\;.
  \end{equation}
\end{lemma}
\begin{proof}
  The mapping $\SP$ is continuous as a superposition of continuous mappings defined by \cref{alg:truncation-projection}.
  To prove the implication denote $\hat \vx \coloneq \SP(\vx)$ and note that $x_i\ge \hat x_i$ for any $i\in [n]$ by construction.
  The latter inequality holds as equality and, therefore, trivially satisfies the implication, unless there exist a constraint $j\in [m]$ such that $\hat x_{n+j}=0$.
  Consider
  \begin{equation}
    \epsilon
    \ge \underbrace{\sum_{i'\in\bar K_j}x_{i'}}_{\le 1+\epsilon}
    -\underbrace{\sum_{i'\in\bar K_j}\hat x_{i'}}_{= 1}
    = \sum_{i'\in K_j}\underbrace{(x_{i'}-\hat x_{i'})}_{\ge 0}
    + \underbrace{x_{n+j}}_{\ge 0}
    - \underbrace{\hat x_{n+j}}_{=0}
    \ge \sum_{i'\in K_j}\vert x_{i'}-\hat x_{i'} \vert
    \stackrel{i\in K_j}{\ge} \vert x_i-\hat x_i \vert \,.
  \end{equation}
\end{proof}
\Crefrange{lem:truncation-projection-is-feasible}{lem:truncation-projection-is-continuous} trivially result in
\begin{theorem}
  Mapping $\SP(\vx)$ defined by \cref{alg:truncation-projection} is a projection to~$\SC$.
\end{theorem}

\section{Sparsity-based speedup}\label{sec:sparsity-based-speedup}

As the temperature gets lower, the majority of values of the primal solution $\vx = \exp(\nicefrac{\vc^\vlambda}{T})$ get closer to zero exponentially fast.
Their influence on the value of the sum $s = \sum_{i\in \bar K_j}x_i$ of all variables in a given clique reduces accordingly, see \cref{alg4:line3} in \cref{alg:exp-domain-stabilized}.
Respectively, most of them can be completely ignored during computation to speed-up the latter.
Below we consider two methods making use of this observation.

\paragraph{Heuristic truncation.}
In our first method we excluded from the summation all variables in a given clique if their values were at least a factor $10^{-8}$ smaller than the largest element in the clique.
This reduced the number of actually processed variables in each clique up to 10 times leading to significant speedup.

However, such a heuristic speedup may lead to an uncontrolled increase of the smoothed dual value $D^T$ ruining algorithms convergence.
This behavior is well-known in the optimal transport domain, see~\cite[Example~3.5]{schmitzer2019stabilized-arxiv}.
Hence, below we derive a method that takes this issue into account.

\paragraph{Truncation based on a bounded dual perturbation (accurate truncation).}
Similarly to~\cite{schmitzer2019stabilized}, our \emph{accurate truncation} is based on the stabilization technique that guarantees a limited maximum value of the exponentiated dual variables $\valpha$, see \cref{alg:exp-domain-stabilized}.
Contrary to~\cite{schmitzer2019stabilized}, we compute the allowed truncation value based on the recent improvement of the smoothed dual value $D^T$ instead of the primal-dual gap.

According to~\cref{lem:smooth-dual-explicit-formula} the dual value $D^T$ is equal to
\begin{equation}
  D^T(\vlambda)=\sum_{j=1}^{m}\lambda_j+T\sum_{i=1}^{n+m}\exp(\nicefrac{c_i}{T})\prod_{j\in J_i}\exp(-\nicefrac{\lambda_j}{T}) = \sum_{j=1}^{m}\lambda_j+T\sum_{i=1}^{n+m}x_i\prod_{j\in J_i}\alpha_j\,.
\end{equation}
Consider the truncated primal estimates
$
  \bar x_i  = x_i\cdot\llbracket x_i \ge \epsilon_i\rrbracket
$
and the respective smoothed dual value
$
  \bar D(\vlambda) \coloneq \sum_{j=1}^{m}\lambda_j+T\sum_{i=1}^{n+m}\bar x_i\prod_{j\in J_i}\alpha_j
$.
Let $\delta$ be a maximum error in computation of $D^T$ allowed due to truncation:
\begin{equation}
  \delta \coloneq D^T(\vlambda)-\bar D(\vlambda)
  = T \sum_{i=1}^{n+m}(x_i-\bar x_i)\prod_{j\in J_i}\alpha_j
  \le T \sum_{i=1}^{n+m}\epsilon_i\prod_{j\in J_i}\tauStab
  = T \sum_{i=1}^{n+m}\epsilon_i\tauStab^{|J_i|}\,.
\end{equation}
Hence, setting
\begin{equation}\label{equ:truncation-threshold}
  \epsilon_i \le \frac{\delta}{T(n+m)\tauStab^{|J_i|}}
\end{equation}
guarantees that the truncation error in computation of the smoothed dual value does not exceed $\delta$.
In our experiments we set $\delta$ to 0.1 of the smoothed dual value improvement attained in the previous computation batch.

\section{Simple primal heuristic}\label{sec:primal-heuristics}

Although the primary subject of our work is the dual optimization, we have implemented a simple primal heuristic to be able to compare to existing methods that return solutions of the non-relaxed MWIS problem.
Our primal heuristic consists of two iteratively repeating steps: \emph{(i)}~A randomized \emph{greedy} feasible solution \emph{proposal generation} and \emph{(ii)}~its \emph{optimized recombination} with the best found solution so far.

\paragraphRunIn{Our greedy proposal generation step} makes use of the dual optimization by being applied to the problem with the reduced costs $\vc^{\vlambda}$ instead of the original ones $\vc$.

\begin{algorithm}[t]
  \addtolength{\hsize}{1.5em}%
  \DontPrintSemicolon
  \KwIn{$\vc^{\vlambda} \in \BR^{n+m}$ \tcc{reduced costs}}
  $\vx \coloneq \mathbf{-1}\in\BR^{n+m}$ \tcc{initialize $\vx$ as undefined}
  \For(\tcc*[h]{iterate over all cliques in random order}){$j \in \mathrm{shuffle}([m])$}{
    \If(\tcc*[h]{proceed if the clique is not set yet}){$\nexists i\in\bar K_j\colon x_i=1$}{
      $i^* \coloneq \argmax_{i \in \bar K_j \text{ s.t. } x_i = -1} c^{\vlambda}_{i}$
      \tcc{find assignable \& locally optimal node}
      $x_{i^*} \coloneq 1$ \tcc{set it to one}
      \For{$ i\colon \{i^*,i\}\in\SE $}{
        $x_i \coloneq 0$
        \tcc{\ldots\ and all its graph neighbors to zero}
      }
    }
  }
  \KwResult{$\vx$}
  \caption{Greedy solution generation.}
  \label{alg:greedy-generation}
\end{algorithm}

The pseudo-code of the generation step is represented by \cref{alg:greedy-generation}.
It considers a graph with additional nodes, one per clique, corresponding to the slack variables introduced in~\eqref{equ:clique-relaxation}.
The algorithm randomly selects a yet unprocessed clique and adds the node with the highest cost in the clique to the solution proposal.
This node and all adjacent nodes are removed from the graph.
The process is repeated until no nodes are left.

\paragraphRunIn{The optimal recombination (also known as crossover or fusion)} step solves the MWIS problem~\eqref{equ:MWIS-edge} only for those variables whose values \emph{differ} in the best integer solution found so far and in the newly generated solution proposal.
The values of all other variables are fixed to the values in best integer solution found so far.
It is known~\cite{eremeev2014optimal} that such optimal recombination problem is reducible in linear time to the min-st-cut/max-flow and is efficiently and exactly solvable therefore.
For implementation reasons we represent the recombination problem as a quadratic pseudo-boolean function and use a dedicated software for its optimization, see \cref{sec:qpbo-transformation} for details.

A similar primal algorithm for the MWIS problem has been proposed, \eg, in~\cite{dong2025metaheuristic}. However, it uses
\emph{(i)}~a different proposal generation, as explained in the following paragraph, and
\emph{(ii)}~multiple local search algorithms applied to each generated solution proposal prior to recombination.

\paragraph{Related work on proposal generation.}
\Cref{alg:greedy-generation}
resembles the WG algorithm of~\cite[Thm.~3]{kako2005approximation} that works with non-negative costs.
Let $\delta_i$ be the node degree.
The WG algorithm iteratively and greedily adds to the solution the node $i^*=\argmax_{i\in\SV}c_i/\delta_i$ and removes it and its neighbors from the graph.
Hence, node degrees change on each step of the algorithm, and one has to maintain a priority queue, that slows down the algorithm.
To speed up computations, this generation method is used to generate only the very first feasible solution proposal in~\cite{dong2025metaheuristic}.
For further solution proposals the changes in the node degrees during generation are not taken into account.

Our proposal generation is much faster than the WG algorithm, since we do not take the node degrees and their respective changes into account.
Importantly, as the dual iterates $\vlambda$ get closer to the dual optimum, our algorithm becomes closer to the WG algorithm.
This is because the optimal cost within each clique gets nearly zero (cf.
\cref{alg:log-sinkhorn}), just like its weighted degree~$c_i/\delta_i$.

\subsection{Discussion: Usage of the relaxed solution}\label{sec:usage-of-relaxed-solution}
Our primal heuristic leverages the relaxation by using reduced costs $\vc^{\vlambda}$ in place of the original costs~$\vc$. In contrast, \cite{dong2025metaheuristic} incorporate the primal relaxed solution $\vx \in [0,1]^n$ into their local search, motivated by the intuition that nodes with high relaxed values are more likely to appear in high-quality integer solutions.

While this is a reasonable heuristic, it lacks theoretical justification in the case of the clique cover relaxation. As follows from \cref{prop:dual-optimality-condition}, all non-zero fractional values in the optimal primal solution correspond to zero reduced costs. Moreover, for duality-based methods such as ours, optimal reduced costs are non-positive -- meaning all non-zero coordinates of $\vx$ correspond to zero reduced costs. Conversely, negative reduced costs correspond to zero-valued primal variables.

This suggests that reduced costs more precisely reflect the dual optimization landscape and are better suited for guiding heuristics. While \cite{dong2025metaheuristic} report empirical gains from using primal relaxed solutions, we suspect that the primary benefit stems from identifying zero and non-zero primal LP values -- information also contained in the reduced costs. Nevertheless, if one prefers to use primal estimates, our truncation projection (\cref{alg:truncation-projection}) offers such estimates even before full dual convergence.

That said, reduced costs correspond to the \emph{modified} objective in~\eqref{equ:clique-relaxation}, so standard MWIS heuristics are not directly applicable and may require non-trivial adaptation, which is not always feasible.

\section{Experimental evaluation}\label{sec:experiments}

\subsection{Datasets}

We focus on large-scale, real-world MWIS problems, in contrast to prior work that primarily evaluated on artificially weighted graphs -- typically unweighted instances with randomized node weights~\cite{lamm2019exactly}. This practice suffers from two key limitations: \emph{(i)}~variability in random weights across studies hinders comparability, and \emph{(ii)}~the synthetic nature of weights limits the relevance of results to real-world applications. Even evaluations on non-random weighted instances~\cite{grossmann2023finding} offer limited insight unless those instances correspond to practical scenarios.
All datasets in our evaluation are readily represented as clique graphs and we reuse these cliques in our algorithms.

\begin{table}
  \usekomafont{caption}%
  \centerline{%
    \begin{tabular}{l c r@{\,}l@{ }r@{ }r@{\,}l r@{\,}l@{ }r@{ }r@{\,}l r@{\,}l@{ }r@{ }r@{\,}l c}
      \toprule
      \bfseries dataset & \bfseries \#\,instances & \multicolumn{5}{c}{\bfseries \#\,nodes} & \multicolumn{5}{c}{\bfseries \#\,edges} & \multicolumn{5}{c}{\bfseries \#\,cliques} & \multicolumn{1}{c}{\bfseries avg.\ clique size}                                                             \\\midrule
      AVR-large  & 17 & 127 & k & -- & 882 & k & 43  & M & -- & 344 & M & 5.4 & k & -- & 38 & k & 120.2  \\
      AVR-medium & 16 & 10  & k & -- & 84  & k & 126 & k & -- & 39  & M & 1.8 & k & -- & 48 & k & 24.4 \\
      AVR-small  & 5  & 979 &   & -- & 14  & k & 2.4 & k & -- & 44  & k & 805 &   & -- & 15 & k & 3.2  \\
      MSCD       & 21 & 5   & k & -- & 8   & k & 66  & k & -- & 198 & k & 22  & k & -- & 36 & k & 8.8  \\
      \bottomrule
    \end{tabular}}
  \caption{
   \textbf{Datasets used for evaluation.}
    We split the \emph{AVR} dataset into \emph{AVR-small}, \emph{AVR-medium} and \emph{AVR-large} based on the number of nodes.
    Here \textbf{\#\,nodes} and \textbf{\#\,edges} stand for the number of nodes and edges in the graph representation of the problem respectively, \textbf{\#\,cliques} gives the number of clique constraints in the respective ILP representation~\eqref{equ:MWIS-clique}, and \textbf{avg.\ clique size} is the average clique size per instance averaged over the number of instances in the dataset.
  }
  \label{tab:datasets}
\end{table}

\emph{Amazon Vehicle Routing~\cite{dong2021new}.}
This dataset contains by far the largest MWIS instances which are publicly available, see \cref{tab:datasets}.
We split this dataset into three subsets, with small~(\emph{AVR-small}), medium~(\emph{AVR-medium}) and large~(\emph{AVR-large}) problem instances, see \cref{tab:datasets} for details.

\emph{Meta-Segmentation for Cell Detection (MSCD)}
This dataset is derived from semi-automated labeling problems for cell segmentations~\cite[Ch.~7]{prakash2022fully} provided by the authors of this work.
We make this real-world MWIS dataset publicly available for usage by other researchers, see \emph{URL not included for peer-review}.
In each instance the nodes represent segmentations that are generated by multiple algorithms and, therefore, overlap with each other.
The MWIS problem consists in selecting their non-overlapping subset with the maximal total cost.
The size of the instances from this dataset lies between \emph{AVR-medium} and \emph{AVR-small}, see \cref{tab:datasets}.

\subsection{Algorithms}
We compare our method to the following solvers:
\begin{itemize}
  \item
  \emph{LP/ILP-Gurobi simplex\,/\,barrier} -- \cite{gurobi} optimizer in single-threaded LP/ILP mode with simplex or barrier method respectively.
    The ILP problem~\eqref{equ:MWIS-clique} and its LP relaxation are considered.
    We use default parameters in all cases, except for the LP modes, where we set the optimality gap to zero to prevent premature termination.
    We selected Gurobi as a representative of a generic LP/ILP solver as it consistently demonstrates strong performance on a wide range of problem instances.
   \item
    \emph{KaMIS branch\,\&\,reduce\,/\,reduce\,\&\,local-search} --
    The framework~\cite{kamis-code} addresses MWIS problems with reduction rules and branching on variables (\emph{branch\,\&\,reduce}), or by running a local search method after an initial reduction (\emph{reduce\,\&\,local-search}).
    We run both with default parameters.
  \item
    \emph{KaMIS local-search only} -- a representative of the local search paradigm for MWIS solvers.
    This variant was obtained by modifying the code to disable the initial reduction step.
  \item
    \emph{METAMIS~\cite{dong2025metaheuristic}, METAMIS+LP~\cite{dong2025metaheuristic}, ILSVND~\cite{nogueira2018hybrid}} -- recent primal MWIS heuristics.
    In contrast to other algorithms, we did not execute these methods ourselves due to the lack of available source code.
    Instead, we compare our results with those presented in~\cite{dong2025metaheuristic}.
\end{itemize}

We note that \cite{grossmann2025concurrent} recently introduced a new heuristic, which surpasses all previously known MWIS approaches, including ours, on the \emph{AVR} dataset described below. Nevertheless, our work remains important as it targets efficient solving of the clique cover LP relaxation, providing strong dual bounds and reduced costs that can complement and enhance such advanced heuristics.

Experiments were conducted on an Intel Core i7-4770 CPU~(3.40\,GHz), measuring wall-clock time.
To ensure reliable and reproducible results, each experiment was run six times and the best result is reported.

\begin{figure}
  \usekomafont{caption}%
  \centerline{%
    \includegraphics{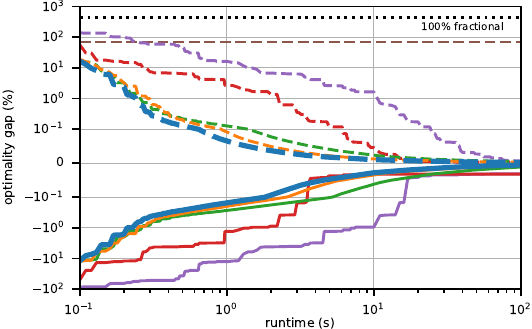}\hfill%
    \includegraphics{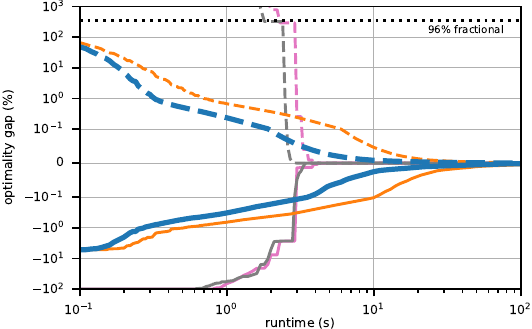}}
  \centerline{
    \makebox[3.54in][c]{\textbf{(a)} AVR-medium (ablation study of our LP-Algorithms)}\hfill%
    \makebox[3.54in][c]{\textbf{(b)}~MSCD}}
  \vspace{1em}
  \centerline{%
    \includegraphics{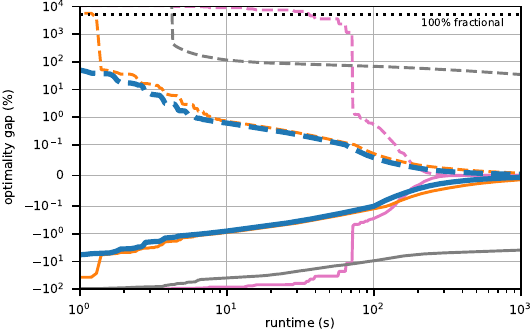}\hfill%
    \includegraphics{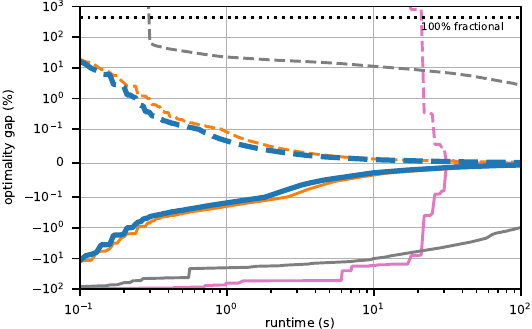}}
  \centerline{%
    \makebox[3.54in][c]{\textbf{(c)}~AVR-large}\hfill%
    \makebox[3.54in][c]{\textbf{(d)}~AVR-medium}}
  \vspace{1em}
  \centerline{%
    \includegraphics{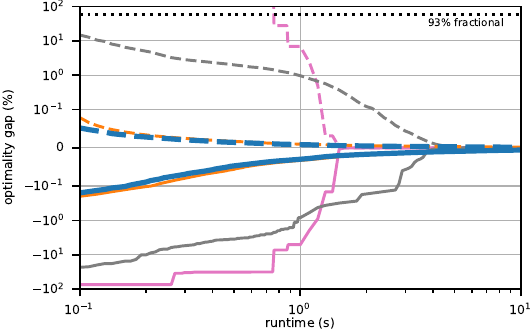}\hfill%
    \includegraphics{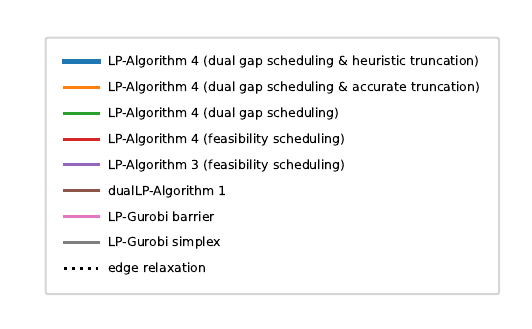}}
  \centerline{%
    \makebox[3.54in][c]{\textbf{(e)}~AVR-small}\hfill}
  \caption{
    \textbf{Optimization of the problem relaxation.}
    These plots track the relative gap to the optimum of the LP relaxation of MWIS problem~\eqref{equ:MWIS-clique} across run time.
    Dashed lines refer to the dual objective (upper bound) and solid lines refer to the primal (relaxed) objective (lower bound).
    Instead of individual problem instances, we show the mean objective across all instances of the respective dataset, providing a consolidated view of algorithm performance.
    Our best algorithms are able to provide high-quality solutions ($<$\,1\,\% rel.~gap) for the majority of the instances within 1--10~seconds. \emph{LP-Gurobi simplex} is run twice, as primal and dual simplex respectively, to obtain both curves.
  }
  \label{fig:plots-relaxed}
\end{figure}

\subsection{Ablation study for the LP relaxation (\cref{fig:plots-relaxed}(a))}
On the \emph{AVR-medium} dataset we compare performance of \cref{alg:non-smooth-bcd,alg:log-sinkhorn,alg:exp-domain-stabilized}
with both \emph{feasibility} and \emph{dual gap scheduling} (\cref{sec:feasibility-constraint-satisfaction,sec:duality-gap-scheduling} resp.)
with \emph{heuristic} or \emph{accurate truncation} (\cref{sec:sparsity-based-speedup}) and without it.

We run all algorithms in LP mode, without the primal heuristic.
The primal bound for them is computed with \cref{alg:truncation-projection}.
All algorithms use the same parameter values, where applicable: $\tauBatch=50$, $\tauDrop=0.5$, $\tauFeas=0.01$, $\tauGap=0.5$, starting temperature $T=0.01$.
All algorithms but those that utilize \emph{accurate speedup} use $\tauStab=10^{30}$, otherwise $\tauStab=10$ to allow for an efficient truncation according to~\eqref{equ:truncation-threshold}.

We report primal-dual convergence plots averaged over all instances in each dataset.
Due to instance similarity, these averaged plots closely resemble individual per-instance plots. Other datasets exhibit similar behavior to \cref{fig:plots-relaxed}(a) and are thus omitted for brevity.

As a baseline, we include the dual bound from the edge relaxation (dotted line), which is consistently loose across \emph{all} datasets (see \cref{fig:plots-relaxed}(b–e)). On average, 95\% of solution values remain fractional (median: 99.4\%), rendering persistency-based graph reduction ineffective.

A second baseline, the non-smoothed dual coordinate descent method (\cref{alg:non-smooth-bcd}), improves upon edge relaxation but remains far inferior to the other algorithms.

As expected, switching from \emph{log-} to \emph{exp-}domain (\ie, from \cref{alg:log-sinkhorn} to \cref{alg:exp-domain-stabilized}) yields nearly an order-of-magnitude speedup. Replacing \emph{feasibility-based scheduling} with \emph{dual gap scheduling} (\cref{sec:T-scheduling}) brings further gains. Feasibility scheduling requires a small threshold $\tauFeas$, as larger values cause slow or stalled convergence due to poor initialization after temperature drops. In contrast, dual gap scheduling dynamically adapts temperature based on solution accuracy.

Finally, both \emph{heuristic} and \emph{accurate truncation} schemes yield additional speedups. Our two most effective configurations -- \cref{alg:exp-domain-stabilized} with duality gap scheduling and either truncation mode -- are used throughout the remaining experiments with the same parameters as above.

\subsection{Our best methods against generic LP solver (\cref{fig:plots-relaxed}(b-e))}
Our next experiment compares performance of \cref{alg:exp-domain-stabilized} with duality gap scheduling and heuristic/accurate truncation against \emph{Gurobi simplex/barrier}.
We do this with the same averaged plots as \cref{fig:plots-relaxed}(a).

The plots show that our methods significantly outperform \emph{Gurobi} at the early optimization stages, but require much longer than \emph{LP-Gurobi barrier} if a very high accuracy is required.
This is, however, a standard behavior of all first order methods including ours, due to their in total sub-linear (more precisely, exponentially low) convergence rate.

Although \cref{alg:exp-domain-stabilized} with \emph{heuristic truncation} is always somewhat faster than with the \emph{accurate truncation}, the difference is negligible in all datasets but \emph{MSCD}.
This is because instances of this dataset require stabilization more often, especially, when the relatively small stabilization threshold $\tauStab=10$ (see~\cref{alg:exp-domain-stabilized}) is used, as this is the case for the \emph{accurate truncation}.
This is the reason why we use the much higher value $\tauStab=10^{30}$ for all other algorithms.

Notably, \emph{LP-Gurobi barrier} significantly outperforms \emph{LP-Gurobi simplex} on all datasets but \emph{MSCD} (for the latter their performance is comparable), hence we stick to \emph{LP-Gurobi barrier} for our later experiments in the ILP mode.

\begin{figure}
  \usekomafont{caption}%
  \centerline{\includegraphics{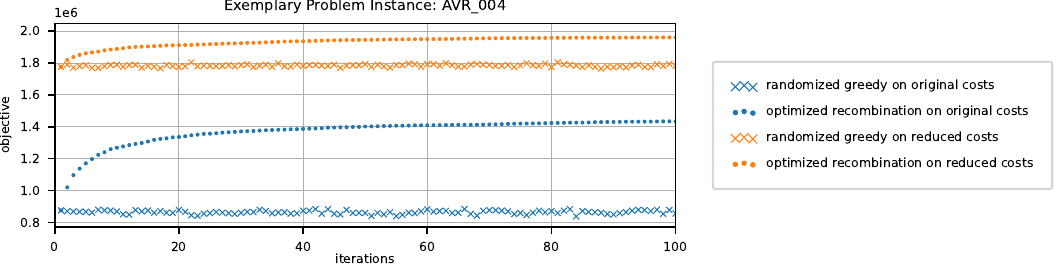}}
  \caption{
    \textbf{Primal heuristic ablation.}
    Exemplary comparison on the problem instance \emph{AVR\_004} from the \emph{AVR-large} dataset.
    As the comparison suggests, optimized recombination leads to significant improvement of the resulting solution.
    Even more improvement is obtained due to usage of the reduced costs:
    Even the best recombination result obtained for the original costs is worse than each single greedily generated solution based on reduced costs.
  }
  \label{fig:plot-fusion}
\end{figure}

\subsection{Ablation study for our primal heuristic (\cref{fig:plot-fusion})}

In our third experiment, we demonstrate the importance of optimized recombination and the use of reduced costs instead of the original ones. This effect is not limited to this example but occurs across most problem instances.

As shown in \cref{fig:plot-fusion}, optimized recombination significantly improves solution quality. An even greater improvement is achieved by using reduced costs: The best recombination result with original costs is inferior to every single greedily generated solution based on reduced costs.

\begin{figure}
  \usekomafont{caption}%
  \centerline{%
    \includegraphics{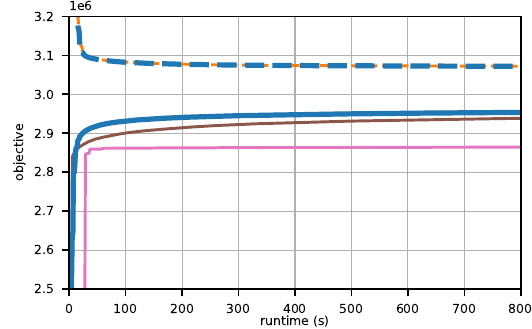}\hfill%
    \includegraphics{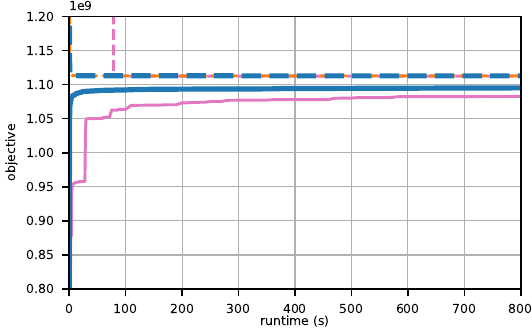}}
  \centerline{%
    \makebox[3.54in][c]{\textbf{(a)}~AVR-large}\hfill%
    \makebox[3.54in][c]{\textbf{(b)}~AVR-medium}}
  \vspace{.5em}
  \centerline{%
    \includegraphics{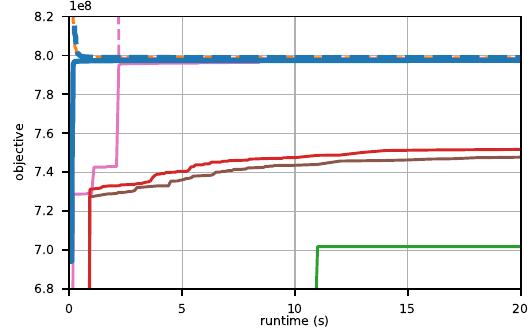}\hfill%
    \includegraphics{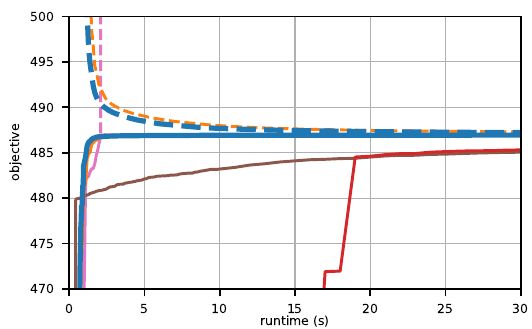}}
  \centerline{%
    \makebox[3.54in][c]{\textbf{(c)}~AVR-small}\hfill%
    \makebox[3.54in][c]{\textbf{(d)}~MSCD}}
  \vspace{1em}%
  \centerline{%
    \includegraphics{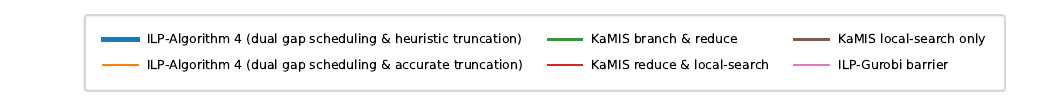}}
  \caption{
    \textbf{Optimization of the non-relaxed problem.}
    These plots present the absolute primal objective values over run time, highlighting the effectiveness of our best-performing dual optimizer in conjunction with the proposed primal heuristics.
    Dashed lines refer to the dual objective (upper bound) and solid lines refer to the primal objective (lower bound).
    We compare our solver to \emph{Gurobi ILP} and \emph{KaMIS}.
    Similarly to \cref{fig:plots-relaxed} we show the mean objective across all problem instances of the different datasets.
    Our solver yields high-quality MWIS solutions after only a few seconds for the majority of instances.
    Several variants of the \emph{KaMIS} algorithms are missing in some plots as their attained values lie outside the visible area of the respective plots.
  }
  \label{fig:plots-integer}
\end{figure}

\subsection{Our best methods against competitors in ILP setting (\cref{fig:plots-integer})}
In this experiment, we compare different methods addressing the \emph{non-relaxed} MWIS problem~\eqref{equ:MWIS-clique}.
As with \cref{fig:plots-relaxed}, we show averaged plots and state their qualitative similarity to individual plots within each dataset.
From our algorithms we show only \cref{alg:exp-domain-stabilized} with duality gap scheduling and heuristic/accurate truncation, where additionally 50~integer solution proposals are generated with \cref{alg:greedy-generation} each $\tauBatch=50$ dual iterations and sequentially recombined with the current best integer solution.

In the short run our algorithms outperforms all its competitors on all datasets but \emph{MSCD}, where the results in terms of the primal integer solution are on par with \emph{ILP-Gurobi barrier}.
In general \emph{ILP-Gurobi barrier} is the most competitive among other methods.
Only \emph{KaMIS local search only} was able to show better results and only for the \emph{AVR-large} dataset.
On all other datasets, each of the three tested KaMIS variants is noncompetitive, at least in the considered time span of 800~seconds.

On the small-sized datasets \emph{AVR-small} and \emph{MSCD} the comparison between our \cref{alg:exp-domain-stabilized} and \emph{ILP-Gurobi barrier} leads to similar results as for the comparison in the LP mode, see \cref{fig:plots-relaxed}.
Specifically, our method is faster initially, quickly finds a good approximate solution, but does not converge to the ILP optimum.
In contrast, \emph{ILP-Gurobi barrier} starts slightly slower but eventually finds the ILP optimum.

However, for larger datasets \emph{AVR-large} and \emph{AVR-medium} \emph{ILP-Gurobi barrier} is not always able to find the optimum and often is even unable to reach the solution accuracy of our methods within 1,000~seconds.

This type of performance is typical when comparing exact methods like \emph{ILP-Gurobi barrier} to approximate ones like ours.
Our method performs well because the reduced costs from dual optimization give it a significant advantage early on.
Being very important for the performance of the primal heuristic at the beginning, its effect diminishes as we approach the solution of the relaxed problem.
In other words, primal heuristics generally cannot fully benefit from high-precision relaxed solutions, whereas Gurobi's branch-and-bound method can.

\subsection{Comparison to METAMIS~\cite{dong2025metaheuristic} and ILSVND~\cite{nogueira2018hybrid}}

\cref{tab:metamis-comparison} shows comparison of our approach against state-of-the-art primal heuristics on all three \emph{AVR} datasets, following the experimental setup of~\cite{dong2025metaheuristic}. Specifically, we run the fastest variant of \emph{Our} method --\cref{alg:exp-domain-stabilized} with our primal heuristic, duality gap scheduling, and heuristic truncation -- five times with different random seeds for 3,600 seconds, matching the time budget used in~\cite{dong2025metaheuristic}. As in their protocol, we report the best result per run for each instance.

\begin{table}
  \fontsize{9pt}{9.5pt}\sffamily
  \newcommand\rot[1]{\rotatebox{90}{#1}}
  \newcommand\Header[1]{\multicolumn{1}{c}{#1}}
  \newcommand\HeaderSeconds{t,~sec}
    \centerline{%
      \begin{tabular}{ll rrrr rrrr}
        \toprule
        && \multicolumn{2}{c}{Our (3,600 sec)}
        & \multicolumn{2}{c}{``Cold'' (\cf \cite[Tab.\,2]{dong2025metaheuristic})}
        & \multicolumn{3}{c}{``Warm'' (\cf \cite[Tab.\,3]{dong2025metaheuristic})}
        \\
        \cmidrule(r){3-4} \cmidrule(lr){5-6} \cmidrule(l){7-9}
        &&& time in \emph{s} until & \multicolumn{2}{c}{time in \emph{s} to surpass \ldots} & \multicolumn{3}{c}{time in \emph{s} to surpass \ldots}
        \\
        \Header{Dataset}
        & \Header{Instance}
        & \Header{Objective}
        & \Header{LP-Gap 0.1\,\%}
        & \Header{ILSVND}
        & \Header{METAMIS}
        & \Header{ILSVND}
        & \Header{METAMIS}
        & \Header{METAMIS+LP}
        \\
        \midrule
        \multirow{17}{*}{\rot{AVR-LARGE}}
         & CR-S-L-1 & \textbf{5,653,748} & 309 & 4 & 63 & 424 & --- & --- \\
         & CR-S-L-2 & \textbf{5,737,640} & 289 & 12 & 198 & 590 & --- & --- \\
         & CR-S-L-4 & \textbf{5,734,265} & 297 & 8 & 113 & 672 & --- & --- \\
         & CR-S-L-6 & \textbf{3,909,318} & 110 & 6 & 34 & 1,840 & --- & --- \\
         & CR-S-L-7 & \textbf{2,011,084} & 19 & 2 & 15 & 507 & --- & --- \\
         & CR-T-C-1 & \textbf{4,723,108} & 50 & 2 & 13 & 122 & --- & --- \\
         & CR-T-C-2 & \textbf{4,945,899} & 122 & 2 & 22 & 166 & --- & --- \\
         & CR-T-D-4 & \textbf{4,881,646} & 170 & 2 & 34 & 583 & --- & --- \\
         & CR-T-D-6 & \textbf{3,013,249} & 75 & 1 & 11 & 290 & --- & --- \\
         & CR-T-D-7 & \textbf{1,455,976} & 10 & $<$ 1 & 3 & --- & --- & --- \\
         & CW-S-L-1 & \textbf{1,648,272} & 190 & 13 & 85 & --- & --- & --- \\
         & CW-S-L-2 & \textbf{1,722,290} & 186 & 25 & 148 & --- & --- & --- \\
         & CW-S-L-4 & \textbf{1,733,085} & 133 & 36 & 463 & --- & --- & --- \\
         & CW-S-L-6 & \textbf{1,164,850} & 35 & 9 & 137 & --- & --- & --- \\
         & CW-S-L-7 & \textbf{589,334} & 10 & 7 & 90 & --- & --- & --- \\
         & CW-T-C-1 & \textbf{1,330,362} & 15 & 6 & 39 & 1,016 & --- & --- \\
         & CW-T-C-2 & \textbf{937,013} & 14 & 10 & 206 & --- & --- & --- \\
        \midrule
        \multirow{16}{*}{\rot{AVR-MEDIUM}}
         & CW-T-D-4 & \textbf{460,198} & 2 & 1 & 6 & 56 & --- & --- \\
         & CW-T-D-6 & \textbf{459,776} & 2 & $<$ 1 & 7 & 28 & --- & --- \\
         & MR-D-03 & \textbf{1,758,074,154} & 1 & $<$ 1 & $<$ 1 & $<$ 1 & 18 & 80 \\
         & MR-D-05 & \textbf{1,790,179,623} & 1 & $<$ 1 & 1 & $<$ 1 & 2 & 8 \\
         & MR-D-FN & \textbf{1,802,268,039} & 1 & $<$ 1 & 1 & $<$ 1 & 3 & 8 \\
         & MR-W-FN & \textbf{5,386,842,781} & $<$ 1 & $<$ 1 & $<$ 1 & $<$ 1 & $<$ 1 & $<$ 1 \\
         & MT-D-200 & \textbf{287,155,108} & $<$ 1 & $<$ 1 & $<$ 1 & $<$ 1 & $<$ 1 & $<$ 1 \\
         & MT-D-FN & \textit{290,866,943} & $<$ 1 & $<$ 1 & 165 & $<$ 1 & $<$ 1 & $<$ 1 \\
         & MT-W-200 & \textbf{384,056,012} & $<$ 1 & $<$ 1 & $<$ 1 & $<$ 1 & $<$ 1 & $<$ 1 \\
         & MT-W-FN & \textbf{390,869,891} & $<$ 1 & $<$ 1 & $<$ 1 & $<$ 1 & $<$ 1 & $<$ 1 \\
         & MW-D-20 & \textbf{526,688,093} & $<$ 1 & $<$ 1 & 1 & $<$ 1 & 1 & 2 \\
         & MW-D-40 & \textbf{538,944,498} & $<$ 1 & $<$ 1 & 2 & $<$ 1 & 1 & 3 \\
         & MW-D-FN & \textbf{545,720,496} & $<$ 1 & $<$ 1 & $<$ 1 & $<$ 1 & 4 & 6 \\
         & MW-W-05 & \textbf{1,329,653,408} & 18 & 1,699 & 1,699 & 1,699 & 1,699 & 1,699 \\
         & MW-W-10 & 1,300,124,841 & 18 & --- & --- & --- & --- & --- \\
         & MW-W-FN & 1,248,553,336 & 17 & --- & --- & --- & --- & --- \\
        \midrule
        \multirow{5}{*}{\rot{AVR-SMALL}}
         & MR-D-01 & \textbf{1,691,358,724} & $<$ 1 & $<$ 1 & 1 & $<$ 1 & $<$ 1 & $<$ 1 \\
         & MT-D-01 & \textit{238,166,485} & $<$ 1 & $<$ 1 & $<$ 1 & $<$ 1 & $<$ 1 & $<$ 1 \\
         & MT-W-01 & \textit{312,121,568} & $<$ 1 & $<$ 1 & $<$ 1 & $<$ 1 & $<$ 1 & $<$ 1 \\
         & MW-D-01 & \textbf{476,361,236} & $<$ 1 & 18 & 141 & 3 & 18 & 18 \\
         & MW-W-01 & \textit{1,270,305,952} & $<$ 1 & $<$ 1 & $<$ 1 & $<$ 1 & $<$ 1 & $<$ 1 \\
        \midrule
        \multirow{3}{*}{\rot{ALL}}
         & & \multicolumn{2}{l}{\bfseries better objective: 32} \\
         & & \multicolumn{2}{l}{\itshape  equal objective:   4} \\
         & & \multicolumn{2}{l}{          worse objective:   2} \\
        \bottomrule
      \end{tabular}}
  \caption{
    \textbf{Comparison of \cref{alg:exp-domain-stabilized} with our primal heuristic, duality gap scheduling and heuristic truncation denoted as \emph{Our} to METAMIS~\cite{dong2025metaheuristic} and ILSVND~\cite{nogueira2018hybrid} based on \cite[Tab.\,4-7]{dong2025metaheuristic}.}
    All algorithms have been run for 3,600 seconds.
    Datasets are sorted by their size, from the largest on top to the smallest in the bottom.
    \textbf{``Warm''} and \textbf{``Cold''} columns correspond to competing algorithms that use or, respectively, do not use a warm-start with a specific (not used by \emph{Our}) initial solution as explained in \cite{dong2025metaheuristic}.
    \textbf{LP-Gap 0.1\%} -- time required to attain the 0.1\% relative duality gap for the relaxed problem. \textbf{ILSVND, METAMIS, METAMIS+LP} -- time required by \emph{Our} algorithm to attain the best value of the respective competing method. This value is taken from the column ``w'' in \cite[Tab.\,4-7]{dong2025metaheuristic}.
    To obtain the full time of \emph{Our} algorithm its time should be added to the respective value in the ``LP-Gap 0.1\,\%'' column. Objective value marked as \textbf{bold} means that this is \emph{strictly greater} than those obtained by any competing algorithm with the ``Cold''-start. \emph{Italic} font means that our value and value of the best competing algorithm coincide: For two of them, MT-D-01 and MT-W-01, this value furthermore coincides with the LP bound and is, therefore, optimal.
    Although \emph{Our} method only uses a comparably simple primal heuristic, it outperforms its competitors due to the usage of the reduced costs, see a detailed explanation in the main text.
  }
  \label{tab:metamis-comparison}
\end{table}

To isolate the contribution of our dual method, we first run dual updates until the relaxed LP reaches a 0.1\% relative duality gap, after which we begin executing the primal heuristic (with 50 proposals per batch).
Starting from that point,
due to sparsity, dual updates typically account for less than 5\% of total runtime and can be largely ignored.

\Cref{tab:metamis-comparison} summarizes the time our method requires to reach or exceed the best results of its competitors. We follow~\cite{dong2025metaheuristic} in evaluating both “cold”- and “warm”-start variants of their algorithms. The latter leverages an initial solution not used by our algorithm. For full objective values and setup details of METAMIS and ILSVND, we refer to their paper.

Our method matches or outperforms both competitors (under cold start) on all but two datasets (MW-W-10 and MW-W-FN). For most \emph{AVR-small} and \emph{AVR-medium} instances, we reach this level of accuracy within 1 second; a few take up to 7 seconds, and only three exceed 100 seconds after attaining the 0.1\% duality gap for the relaxed problem. On \emph{AVR-large}, our method outperforms METAMIS within 100 seconds on 11 of 17 instances, and never exceeds 206 seconds. In contrast, METAMIS and ILSVND require the full 3,600 seconds for ``larger instances'', see~\cite{dong2025metaheuristic}. Our method also reaches ILSVND's best reported values \emph{significantly} faster than METAMIS, \cf~\cite{dong2025metaheuristic}[Tab.~2, column METAMIS-$t^*$[s]].

Remarkably, \emph{Our} method outperforms even the warm-start versions of METAMIS and ILSVND on most \emph{AVR-small} and \emph{AVR-medium} problems -- despite not using the warm-start solution that significantly boosts their performance. METAMIS+LP further exploits the LP solution as a local search guide, unlike our purely dual-informed heuristic.

\emph{Why does our method outperform METAMIS, despite its stronger heuristic?}
Our experiments suggest that the use of reduced costs, rather than original ones, is key. This effect is clearly visible in \cref{fig:plot-fusion}. However, as noted in~\cref{sec:usage-of-relaxed-solution}, applying reduced costs in MWIS requires adapting the heuristic, and not all heuristics can accommodate this. This prevents us from re-running KaMIS or ILSVND under our reduced-cost regime to further validate this finding.

\subsection{Greedily generated cliques vs.~cliques from the AVR and MSCD datasets}

As noted in \cref{sec:clique-cover-family}, a clique cover can be generated directly from the problem graph, \eg, with the greedy algorithm described therein. Since much of the literature on MWIS assumes that no clique cover is provided as part of the input, we repeated all our experiments with the best-performing variant of \emph{Our} method (\cref{alg:exp-domain-stabilized} with our primal heuristic, duality gap scheduling, and heuristic truncation) using such greedily generated cliques in place of those shipped with the datasets. For the \emph{AVR} dataset, the relaxation based on the generated clique cover is only slightly less tight than the one based on the original cliques. Consequently, this incurs only a minor, almost negligible performance drop of our method. For the \emph{MSCD} dataset, the generated cover yields a looser relaxation, but its description is far more compact than the original one. This, in turn, speeds up our dual solver as well as the primal heuristic and in total leads to a comparably good integer solution even faster than with the original clique cover.
We provide the detailed comparison in \cref{app:generated_vs_original_clique_cover_comparison}.

\section{Conclusion}

We introduced a new method for solving the maximum-weight independent set problem, centered on an efficient approximate solver for the clique cover LP relaxation. By leveraging reduced costs, we simplify and accelerate the primal heuristic, yielding a competitive approximate solver for the original (non-relaxed) problem. Our approach shows particular strength on large-scale instances and is expected to perform even better when combined with more powerful heuristics that incorporate local search.
This scalability stems from the first-order nature of our method; consequently, it remains an open question whether other first-order techniques (see, \eg, \cite{beck2017first}) could be adapted to this setting and achieve comparable performance.

\section{Acknowledgments}
We thank Prof.~Bernhard Schmitzer for an insightful discussion about modern advances in computational optimal transport.
The discussion lead to the decisive speedup of our algorithms described in this work.
A special thanks to Prof.~Paul Swoboda for thorough proofreading and numerous comments that notably improved this paper.
This work was supported by the German Research Foundation projects 498181230 and 539435352.
Authors further acknowledge facilities for high-throughput calculations bwHPC of the state of Baden-Württemberg (DFG grant INST 35/1597-1 FUGG) as well as the Center for Information Services and High Performance Computing (ZIH) at TU Dresden.

\printbibliography

@string{Conf_NeurIPS = "Advances in Neural Information Processing Systems"}

@inproceedings{grossmann2025concurrent,
  title={Concurrent Iterated Local Search for the Maximum Weight Independent Set Problem},
  author={Gro{\ss}mann, Ernestine and Langedal, Kenneth and Schulz, Christian},
  booktitle={23rd International Symposium on Experimental Algorithms (SEA 2025)},
  pages={22--1},
  year={2025},
  organization={Schloss Dagstuhl--Leibniz-Zentrum f{\"u}r Informatik}
}

@misc{kamis-code,
  author={{K}a{MIS}},
  year={2024},
  title={{K}a{MIS} - {K}arlsruhe Maximum Independent Sets
},
  note={\url{https://karlsruhemis.github.io/}}
}

@book{fang2012entropy,
  title={Entropy optimization and mathematical programming},
  author={Fang, Shu-Cherng and Rajasekera, Jay R. and Tsao, H.-S. Jacob},
  volume={8},
  year={2012},
  publisher={Springer Science \& Business Media}
}

@misc{qpbo-code,
  author = {Kolmogorov, Vladimir},
  title={Source code of algorithms for minimizing functions of binary variables with unary and pairwise terms based on roof duality},
  note={\url{https://pub.ista.ac.at/~vnk/software/QPBO-v1.4.src.zip}},
  year={2014}
}

@inproceedings{grossmann2023finding,
  title={Finding near-optimal weight independent sets at scale},
  author={Gro{\ss}mann, Ernestine and Lamm, Sebastian and Schulz, Christian and Strash, Darren},
  booktitle={Proceedings of the Genetic and Evolutionary Computation Conference},
  pages={293--302},
  year={2023}
}

@article{moon1965cliques,
  title={On cliques in graphs},
  author={Moon, John W. and Moser, Leo},
  journal={Israel journal of Mathematics},
  volume={3},
  pages={23--28},
  year={1965},
  publisher={Springer}
}

@book{schrijver2003combinatorial,
  title={Combinatorial optimization: polyhedra and efficiency},
  author={Schrijver, Alexander},
  volume={24},
  number={2},
  year={2003},
  publisher={Springer}
}

@article{pullan2009optimisation,
  title={Optimisation of unweighted\allowbreak/\allowbreak weighted maximum independent sets and minimum vertex covers},
  author={Pullan, Wayne},
  journal={Discrete Optimization},
  volume={6},
  number={2},
  pages={214--219},
  year={2009},
  publisher={Elsevier}
}

@article{paschos1997survey,
  title={A survey of approximately optimal solutions to some covering and packing problems},
  author={Paschos, Vangelis T.},
  journal={ACM Computing Surveys (CSUR)},
  volume={29},
  number={2},
  pages={171--209},
  year={1997},
  publisher={ACM New York, NY, USA}
}

@article{nogueira2018hybrid,
  title={A hybrid iterated local search heuristic for the maximum weight independent set problem},
  author={Nogueira, Bruno and Pinheiro, Rian G.S. and Subramanian, Anand},
  journal={Optimization Letters},
  volume={12},
  pages={567--583},
  year={2018},
  publisher={Springer}
}

@article{shah2005max,
  title={Max product for max-weight independent set and matching},
  author={Shah, Devavrat},
  journal={arXiv preprint cs/0508097},
  year={2005}
}

@article{sanghavi2007message,
  title={Message passing for max-weight independent set},
  author={Sanghavi, Sujay and Shah, Devavrat and Willsky, Alan},
  journal={Advances in Neural Information Processing Systems},
  volume={20},
  year={2007}
}

@article{peyre2017computational,
  title={Computational optimal transport},
  author={Peyr{\'e}, Gabriel and Cuturi, Marco and others},
  journal={Center for Research in Economics and Statistics Working Papers},
  number={2017-86},
  year={2017}
}

@inproceedings{lamm2019exactly,
  title={Exactly solving the maximum weight independent set problem on large real-world graphs},
  author={Lamm, Sebastian and Schulz, Christian and Strash, Darren and Williger, Robert and Zhang, Huashuo},
  booktitle={2019 Proceedings of the Twenty-First Workshop on Algorithm Engineering and Experiments (ALENEX)},
  pages={144--158},
  year={2019},
  organization={SIAM}
}

@article{dong2025metaheuristic,
  title={A metaheuristic algorithm for large maximum weight independent set problems},
  author={Dong, Yuanyuan and Goldberg, Andrew V and Noe, Alexander and Parotsidis, Nikos and Resende, Mauricio GC and Spaen, Quico},
  journal={Networks},
  volume={85},
  number={1},
  pages={91--112},
  year={2025},
  publisher={Wiley Online Library}
}

@article{trevisan2014inapproximability,
  title={Inapproximability of combinatorial optimization problems},
  author={Trevisan, Luca},
  journal={Paradigms of Combinatorial Optimization: Problems and New Approaches},
  pages={381--434},
  year={2014},
  publisher={Wiley Online Library}
}

@inproceedings{kako2005approximation,
  title={Approximation algorithms for the weighted independent set problem},
  author={Kako, Akihisa and Ono, Takao and Hirata, Tomio and Halld{\'o}rsson, Magn{\'u}s M.},
  booktitle={Graph-Theoretic Concepts in Computer Science: 31st International Workshop, WG 2005, Metz, France, June 23-25, 2005, Revised Selected Papers 31},
  pages={341--350},
  year={2005},
  organization={Springer}
}

@incollection{grotschel1984polynomial,
  title={Polynomial algorithms for perfect graphs},
  author={Gr{\"o}tschel, Martin and Lov{\'a}sz, L{\'a}szl{\'o} and Schrijver, Alexander},
  booktitle={North-Holland mathematics studies},
  volume={88},
  pages={325--356},
  year={1984},
  publisher={Elsevier}
}

@article{boros2002pseudo,
  title={Pseudo-boolean optimization},
  author={Boros, Endre and Hammer, Peter L.},
  journal={Discrete applied mathematics},
  volume={123},
  number={1-3},
  pages={155--225},
  year={2002},
  publisher={Elsevier}
}

@article{rodriguez2019persistency,
  title={Persistency of Linear Programming Formulations for the Stable Set Problem},
  author={Rodr{\'i}guez-Heck, Elisabeth and Stickler, Karl and Walter, Matthias and Weltge, Stefan},
  journal={Mathematical programming},
  volume={192},
  issue={1},
  pages={387-407},
  year={2022}
}

@article{nemhauser1975vertex,
  title={Vertex packings: {S}tructural properties and algorithms},
  author={Nemhauser, George L. and Trotter Jr., Leslie E.},
  journal={Mathematical Programming},
  volume={8},
  number={1},
  pages={232--248},
  year={1975},
  publisher={Springer}
}

@incollection{hochbaum1996approximating,
  title={Approximating covering and packing problems: set cover, vertex cover, independent set, and related problems},
  author={Hochbaum, Dorit},
  booktitle={Approximation algorithms for NP-hard problems},
  pages={94--143},
  year={1996}
}

@article{eremeev2014optimal,
  title={Optimal recombination in genetic algorithms for combinatorial optimization problems: {P}art {I}},
  author={Eremeev, Anton and Kovalenko, Julia},
  journal={Yugoslav Journal of Operations Research},
  volume={24},
  number={1},
  pages={1--20},
  year={2014}
}

@article{prusa2019solving,
  title={Solving {LP} relaxations of some {NP}-hard problems is as hard as solving any linear program},
  author={Pr{\r{u}}{\v{s}}a, Daniel and Werner, Tom\'a\v{s}},
  journal={SIAM Journal on Optimization},
  volume={29},
  number={3},
  pages={1745--1771},
  year={2019},
  publisher={SIAM}
}

@inproceedings{haller2022comparative,
  title={A comparative study of graph matching algorithms in computer vision},
  author={Haller, Stefan and Feineis, Lorenz and Hutschenreiter, Lisa and Bernard, Florian and Rother, Carsten and Kainm{\"u}ller, Dagmar and Swoboda, Paul and Savchynskyy, Bogdan},
  booktitle={European Conference on Computer Vision},
  pages={636--653},
  year={2022},
  organization={Springer}
}

@article{dlask2025relative,
  title={Relative-interior solution for the (incomplete) linear assignment problem with applications to the quadratic assignment problem},
  author={Dlask, Tom{\'a}{\v{s}} and Savchynskyy, Bogdan},
  journal={Annals of Mathematics and Artificial Intelligence},
  pages={1--44},
  year={2025},
  publisher={Springer}
}

@article{schmitzer2019stabilized,
  title={Stabilized sparse scaling algorithms for entropy regularized transport problems},
  author={Schmitzer, Bernhard},
  journal={SIAM Journal on Scientific Computing},
  volume={41},
  number={3},
  pages={A1443--A1481},
  year={2019},
  publisher={SIAM}
}

@article{schmitzer2019stabilized-arxiv,
  title={Stabilized sparse scaling algorithms for entropy regularized transport problems},
  author={Schmitzer, Bernhard},
  journal={arXiv Preprint arXiv:1610.06519v1},
  year={2016}
}

@article{kushinsky2019sinkhorn,
  title={Sinkhorn algorithm for lifted assignment problems},
  author={Kushinsky, Yam and Maron, Haggai and Dym, Nadav and Lipman, Yaron},
  journal={SIAM Journal on Imaging Sciences},
  volume={12},
  number={2},
  pages={716--735},
  year={2019},
  publisher={SIAM}
}

@inproceedings{weed2018explicit,
  title={An explicit analysis of the entropic penalty in linear programming},
  author={Weed, Jonathan},
  booktitle={Conference On Learning Theory},
  pages={1841--1855},
  year={2018},
  organization={PMLR}
}

@inproceedings{brendel2010segmentation,
  title={Segmentation as maximum-weight independent set},
  author={Brendel, William and Todorovic, Sinisa},
  booktitle=Conf_NeurIPS,
  year={2010}
}

@phdthesis{prakash2022fully,
  title={Fully Unsupervised Image Denoising, Diversity Denoising and Image Segmentation with Limited Annotations},
  author={Prakash, Mangal},
  year={2022},
  school={Technische Universit{\"a}t Dresden}
}

@inproceedings{dong2021new,
  title={New instances for maximum weight independent set from a vehicle routing application},
  author={Dong, Yuanyuan and Goldberg, Andrew V. and Noe, Alexander and Parotsidis, Nikos and Resende, Mauricio G.C. and Spaen, Quico},
  booktitle={Operations Research Forum},
  volume={2},
  pages={1--6},
  year={2021},
  organization={Springer}
}

@inproceedings{karp2010reducibility,
  title={Reducibility among combinatorial problems},
  author={Karp, Richard M.},
  year={1972},
  publisher={Springer},
  booktitle={Complexity of computer computations}
}

@book{garey1979computers,
  title={Computers and intractability},
  author={Garey, Michael R. and Johnson, David S.},
  year={1979},
  publisher={W.~H.~Freeman and Company, New York}
}

@book{bertsekas1999nonlinear,
  title={Nonlinear programming, second edition},
  author={Bertsekas, Dimitri P.},
  year={1999},
  publisher={Athena scientific}
}

@article{bregman1967relaxation,
  title={The Relaxation Method of Finding the Common Point of Convex Sets and its Application to the Solution of Problems in Convex Programming},
  author={Bregman, Lev M.},
  journal={USSR Computational Mathematics and Mathematical Physics},
  year={1967}
}

@article{sinkhorn1964relationship,
  title={A Relationship Between Arbitrary Positive Matrices and Doubly Stochastic Matrices},
  author={Sinkhorn, Richard},
  journal={The Annals of Mathematical Statistics},
  year={1964}
}

@article{savchynskyy2019discrete,
  author={Savchynskyy, Bogdan},
  title={Discrete Graphical Models -- An Optimization Perspective},
  journal={Foundations and Trends in Computer Graphics and Vision},
  year={2019},
}

@misc{gurobi,
  author={Gurobi},
  title={Gurobi Optimization},
  year={2018},
  note={\url{http://www.gurobi.com}}
}

@article{bron1973algorithm,
  title={Algorithm 457: {F}inding all cliques of an undirected graph},
  author={Bron, Coen and Kerbosch, Joep},
  journal={Communications of the ACM},
  volume={16},
  number={9},
  pages={575--577},
  year={1973},
  publisher={ACM New York, NY, USA}
}

@book{beck2017first,
  title={First-order methods in optimization},
  author={Beck, Amir},
  year={2017},
  publisher={SIAM}
}

\clearpage
\appendix
\crefalias{section}{appendix}

\section{Proof of~\cref{prop:edge-LP-solution-for-complete-graph}}
\label{sec:proof:edge-LP-solution-for-complete-graph}

\edef\oldproposition{\arabic{proposition}}
\setcounter{proposition}{0}
\begin{proposition}
  Let $G$ be fully-connected, all costs positive, $c_j > 0$, $j\in\SV$, and $c_i < \sum_{j\in\SV\backslash\{i\}}c_j$ for $i=\argmax_{j\in \SV}c_j$.
  Then the edge relaxation has a unique solution $(\underbrace{\nicefrac{1}{2},\nicefrac{1}{2},\dots,\nicefrac{1}{2}}_{n})$.
\end{proposition}
\setcounter{proposition}{\oldproposition}
\begin{proof}
  Due to half-integrality, it is sufficient to consider only half-integer feasible solutions.
  Due to the cost non-negativity the feasible solution $\vy=\underbrace{(\frac{1}{2},\frac{1}{2},\dots,\frac{1}{2})}_{n}$ is the best among those with coordinates $1/2$ and $0$ only, as it contains the maximal number of non-zero coordinates.
  The cost of this solution is $\langle \vc,\vy\rangle=\frac{1}{2}\sum_{j\in\SV}c_j$.

  Consider now feasible solutions that are at least one coordinate equal to $1$.
  Since the graph is fully-connected, there are only $|\SV|$ of them, \eg, those of the form 
  $\vx^i=(\underbrace{0,\dots,0}_{i-1},1,0,\dots,0)$ 
  with the total costs equal to $c_i$, respectively.
  The condition of the proposition implies
  \begin{equation}
    c_i < \sum_{j\in\SV\backslash\{i\}}c_j = \sum_{j\in\SV}c_j - c_i
  \end{equation}
  and, therefore, $c_i < \frac{1}{2}\sum_{j\in\SV}c_j = \langle \vc,\vy\rangle$.
\end{proof}

\section{Proof that maximal clique relaxation is also tight for complete graphs}
\label{sec:prop:clique-LP-solution-for-complete-graph}

\begin{proposition}\label{prop:clique-LP-solution-for-complete-graph}
  The maximal clique relaxation, defined by the single maximal clique, is tight for complete graphs.
\end{proposition}
\begin{proof}
  It is sufficient to prove that the feasible set being a polytope has only integer vertices.
  Indeed, assume $\vx'\in [0,1]^n$ is a vertex.
  Then there is $\vc\in\BR^n$ such that $\vx'$ is a unique solution of $\{\max_{\vx\in [0,1]^n}\langle \vc,\vx\rangle, \text{ s.t. } \sum_{i=1}^n x_i \le 1\}$.
  Assume that $\vx'$ is fractional and w.l.o.g.
  $c_1\ge c_2>0$ and $x'_1,x'_2>0$. Then
  \begin{equation}\label{eq:proof-integrality-for-complete-graph}
    \langle \vc,\vx'\rangle = c_1x'_1+c_2x'_2+\sum_{i=3}^n c_ix'_i
    \le c_1(x'_1+x'_2)+c_2\cdot0+\sum_{i=3}^n c_ix'_i = \langle \vc,\vx''\rangle\,,
  \end{equation}
  where $\vx''$ such that $x''_1=x'_1+x'_2$, $x''_2=0$ and $x''_i=x'_i$ for $i=3\dots n$. $\vx''$ is feasible since $\vx'$ is feasible as vertex.
  If the inequality in~\eqref{eq:proof-integrality-for-complete-graph} holds strictly, then $\vx'$ is not a solution and if it holds as equality, then $\vx'$ is not a unique one.
  So we obtained a contradiction.
\end{proof}
\section{Example for loose edge relaxation}
\label{sec:ex:bcd-fix-point}

\begin{example}\label{ex:bcd-fix-point}
  Consider an MWIS problem in format~\eqref{equ:clique-relaxation} with the following $m=6$ constraints:
  \begin{align}
     & x_1+x_2+\textcolor{gray}{x_6} =1     \\
     & x_2+x_3+\textcolor{gray}{x_7} =1     \\
     & x_3+x_1+\textcolor{gray}{x_8} =1     \\
     & x_1+x_3+x_4+\textcolor{gray}{x_9} =1 \\
     & x_4+x_5+\textcolor{gray}{x_{10}} =1  \\
     & x_5+x_3+\textcolor{gray}{x_{11}} =1
  \end{align}
  The coordinates with indices 6 and above correspond to the slack variables and can be found in a single equation only.
  We also assume $c^\vlambda_i=0$ for $i\in[5]$ and $c^\vlambda_i < 0$ for $i > 5$.

  Obviously this $\vlambda$ is a fix-point of \cref{alg:non-smooth-bcd}.
  We will show that this point is not optimal.
  For this, it is sufficient to show that the set $\Sign(\vc^{\vlambda})$ does not contain feasible vectors.
  Assume it does and $\vx\in\Sign(\vc^{\vlambda})$, then $x_i=0$ for $i>5$ as $c^{\vlambda}_i < 0$.
  Assume $x_1=\gamma\ge 0$.
  This implies $x_2=1-\gamma$, and $x_3=\gamma=1-\gamma =0.5$ and hence $x_1=x_2=x_3=0.5$ from the first three constraints.
  The next two constraints imply $x_4=0$ and $x_5=1$.
  The last equality does not hold as $x_5+x_3=1.5 > 1$, which is a contradiction.
\end{example}

\section{Proof of \cref{prop:log-sinkhorn-equiv-to-naive-sinkhorn}}
\label{sec:prop:log-sinkhorn-equiv-to-naive-sinkhorn}

\edef\oldlemma{\arabic{lemma}}
\setcounter{lemma}{1}
\begin{lemma}
  \cref{alg:log-sinkhorn} is equivalent to \cref{alg:naive-sinkhorn} through the transformation $\vx=\exp(\nicefrac{\vc^{\vlambda}}{T})$.
\end{lemma}
\setcounter{lemma}{\oldlemma}
\begin{proof}
  Indeed, let $\lambda'_k=\lambda_k$ for $k\neq j$ and
  \begin{equation}
    \lambda'_j \coloneq \lambda_j+c^{\vlambda}_{i^*}+T\log\sum_{i\in \bar K_j}\exp{\left(\nicefrac{(c^{\vlambda}_i-c^{\vlambda}_{i^*})}{T}\right)}
    = \lambda_j+T\log\sum_{i\in \bar K_j}\exp(\nicefrac{c^{\vlambda}_i}{T})\,.
  \end{equation}
  Taking in account the transformation $\vx=\exp(\vc^{\vlambda}/T)$, the elementary update in \cref{alg:naive-sinkhorn} can be rewritten as
  \begin{align}
    x'_i &\coloneq \frac{x_i}{\sum_{i\in \bar K_j}x_i}
    = \frac{\exp(\nicefrac{c^\vlambda_i}{T})}{\sum_{i\in \bar K_j}\exp(\nicefrac{c^\vlambda_i}{T})}
    = \exp\left(\frac{c^\vlambda_i - T\log \sum_{i\in \bar K_j}\exp(\nicefrac{c^\vlambda_i}{T})}{T} \right) \\
    &= \exp\left(\frac{c_i - \sum_{k\in J_i}\lambda_k - T\log \sum_{i\in \bar K_j}\exp(\nicefrac{c^\vlambda_i}{T})}{T} \right)\\
    &= \exp\left(\frac{c_i - \sum_{k\in J_i\backslash\{j\}}\lambda_k - (\lambda_j + T\log \sum_{i\in \bar K_j}\exp(\nicefrac{c^\vlambda_i}{T}))}{T} \right)\\
    &= \exp\left(\frac{c_i - \sum_{k\in J_i}\lambda'_k}{T} \right)
    = \exp\left(\nicefrac{c^{\vlambda'}_i}{T} \right)\,.
  \end{align}
\end{proof}

\section{Transformation of the recombination problem (\cref{sec:primal-heuristics}) into quadratic pseudo-boolean maximization}
\label{sec:qpbo-transformation}

Let $\SV'$ stand for the set of graph nodes that differ in two feasible solutions $\vx^1, \vx^2\in \{0,1\}^{n}$ of the MWIS problem~\eqref{equ:MWIS-edge} to be fused, \ie, $x^1_i=\bar x^2_i$, $i\in\SV'$ and $x^1_i=x^2_i$ for $i\in\SV\backslash\SV'$ respectively.
Let $\vc$ be the cost vector of the problem.
Let also the set of edges $\SE'$ be inherited from the initial problem: $\SE'=\{\{i,j\}\in\SE\colon i,j\in\SV'\}$.
By construction, the graph $(\SV',\SE')$ is bipartite, that is $\SE'\subseteq\SV^1\times\SV^2$, where $\SV'=\SV^1\cup\SV^2$, $\SV^1\cap\SV^2=\emptyset$, $\SV^k=\{i\in\SV'\colon x^k_i=1\}$, $k=1,2$.
For the sake of notation we assume the edges to be oriented, \eg, $(i,i')\in\SE'$ implies $i\in \SV^1$ and $i'\in\SV^2$.

Consider the variable substitution
\begin{equation}
  \forall i\in \SV'\colon (y_i,w_i) \coloneq
  \begin{cases}
    (x_i, c_i),      & i\in \SV^1    \\
    (\bar x_i,-c_i), & i\in \SV^2\,.
  \end{cases}
\end{equation}
Then for any $i\in\SV^2$ and binary $x_i$ it holds $c_ix_i=w_i y_i - w_i$.
This turns the respective quadratic pseudo-boolean function (see \Cref{sec:relation-to-QPBO}) of the considered MWIS problem into the so called \emph{oriented form}~\cite[Ch.~11]{savchynskyy2019discrete} as follows:
\begin{equation}\label{equ:recombination-oriented-form}
  \sum_{i\in\SV'}c_i x_i - M\cdot\sum_{(i,j)\in\SE'}x_i x_j
  = -\sum_{l\in\SV^2}w_l + \sum_{i\in\SV'}w_i y_i - M\cdot\sum_{(i,j)\in\SE'}y_i \bar y_j\,.
\end{equation}
As shown in, \eg, in~\cite[Ch.~11]{savchynskyy2019discrete}, maximization of the right-hand-side of~\eqref{equ:recombination-oriented-form} reduces to the min-st-cut and we used the software~\cite{qpbo-code} that implements this reduction and solves the respective min-st-cut problem.

\section{Generated vs. original clique cover comparison}\label{app:generated_vs_original_clique_cover_comparison}

\begin{figure}
  \usekomafont{caption}%
  \centerline{%
    \includegraphics{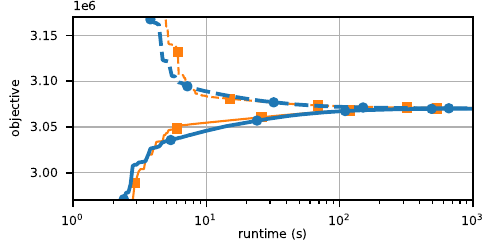}\hfill%
    \includegraphics{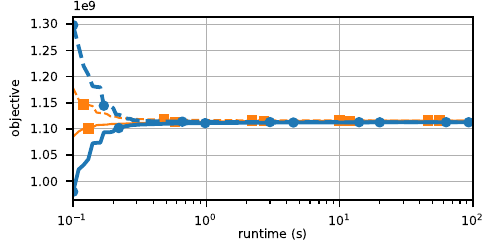}}
  \centerline{%
    \makebox[3.25in][c]{\textbf{(a)}~AVR-large}\hfill%
    \makebox[3.25in][c]{\textbf{(b)}~AVR-medium}}
  \vspace{1em}
  \centerline{%
    \includegraphics{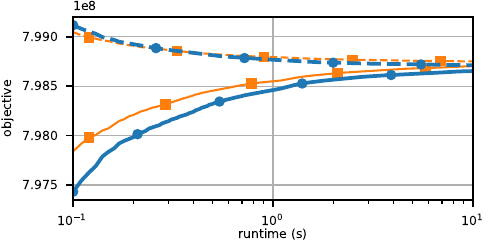}\hfill%
    \includegraphics{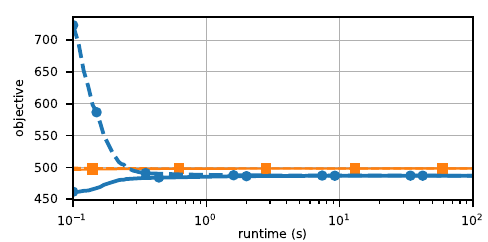}}
  \centerline{%
    \makebox[3.25in][c]{\textbf{(c)}~AVR-small}\hfill%
    \makebox[3.25in][c]{\textbf{(d)}~MSCD}}
  \vspace{1em}
  \centerline{%
    \includegraphics{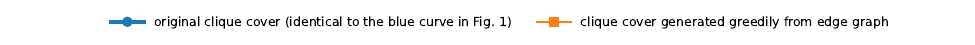}}
  \caption{
    \textbf{Optimization of the problem relaxation with the original and the greedily generated clique cover.}
    These plots track the primal and dual objective values of the clique cover LP relaxation of the MWIS problem~\eqref{equ:MWIS-clique} across run time, obtained by \emph{Our} method (\cref{alg:exp-domain-stabilized} with duality gap scheduling and heuristic truncation), run in LP mode without the primal heuristic.
    Dashed lines refer to the dual objective (upper bound) and solid lines to the primal (relaxed) objective (lower bound).
    Blue curves use the \emph{original} clique cover shipped with the dataset (identical to the corresponding curves in \cref{fig:plots-relaxed}), while orange curves use the clique cover \emph{generated greedily} from the edge graph (\cref{sec:clique-cover-family}).
    Instead of individual instances, we show the mean objective across all problem instances of the respective dataset.
    In most cases the relaxation based on the generated clique cover is weaker than the one based on the original cliques, though the difference is usually minor; for \emph{AVR-medium} and, especially, \emph{MSCD}, the generated cover admits a more compact description, which leads to faster convergence.
    Note that, in contrast to \cref{fig:plots-relaxed}, the vertical axis reports the objective value rather than the relative duality gap. This change  makes the tightness of the two relaxations directly comparable, as their relaxed optima differ.
  }
  \label{fig:plots-relaxed-appendix}
\end{figure}

\begin{figure}
  \usekomafont{caption}%
  \centerline{%
    \includegraphics{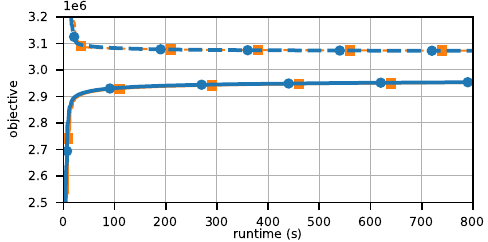}\hfill%
    \includegraphics{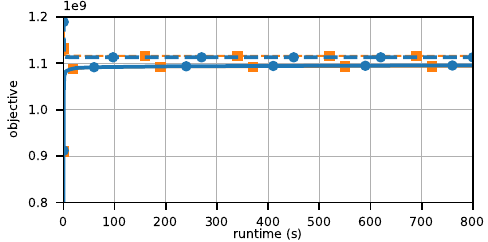}}
  \centerline{%
    \makebox[3.25in][c]{\textbf{(a)}~AVR-large}\hfill%
    \makebox[3.25in][c]{\textbf{(b)}~AVR-medium}}
  \vspace{.5em}
  \centerline{%
    \includegraphics{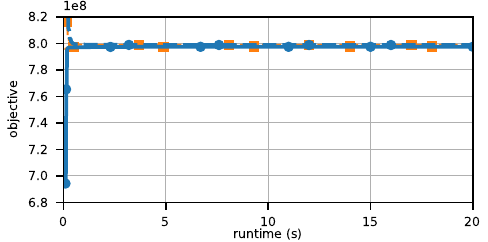}\hfill%
    \includegraphics{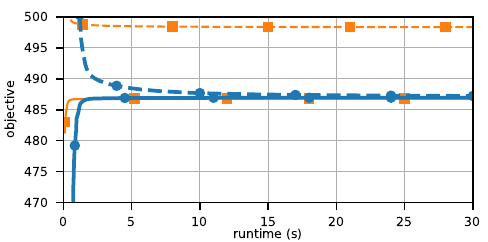}}
  \centerline{%
    \makebox[3.25in][c]{\textbf{(c)}~AVR-small}\hfill%
    \makebox[3.25in][c]{\textbf{(d)}~MSCD}}
  \vspace{1em}%
  \centerline{%
    \includegraphics{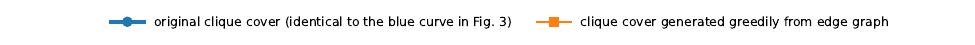}}
  \caption{
    \textbf{Optimization of the non-relaxed problem with the original and the greedily generated clique cover.}
    These plots present the objective values over run time, obtained by \emph{Our} method -- \cref{alg:exp-domain-stabilized} with our primal heuristic, duality gap scheduling, and heuristic truncation.
    Dashed lines refer to the dual objective (upper bound) and solid lines to the primal integer objective (lower bound).
    Similarly to \cref{fig:plots-integer}, we show the mean objective across all problem instances of the respective dataset. Blue curves use the \emph{original} clique cover shipped with the dataset (identical to the corresponding curves in \cref{fig:plots-integer}), while orange curves use the clique cover \emph{sampled greedily} from the edge graph (\cref{sec:clique-cover-family}).
    On the \emph{AVR} datasets the two clique covers yield nearly indistinguishable curves, whereas on \emph{MSCD} the generated cover produces a weaker (looser) dual bound but attains a comparable primal integer solution faster.
  }
  \label{fig:plots-integer-appendix}
\end{figure}

In this section we compare the performance of \emph{Our} method (\cref{alg:exp-domain-stabilized} with our primal heuristic, duality gap scheduling, and heuristic truncation) using greedily generated clique covers against the clique covers shipped with the datasets.
\Cref{fig:plots-relaxed-appendix} reports the comparison for the LP version of the solver (without the primal heuristic), \cref{fig:plots-integer-appendix} for the ILP version, and \cref{tab:metamis-clique-comparison} the relative performance under both clique covers, compared to the METAMIS~\cite{dong2025metaheuristic} algorithm.

\clearpage

\clearpage
\begin{table}
   \fontsize{9pt}{9.5pt}\sffamily
    \newcommand\rot[1]{\rotatebox{90}{#1}}
    \newcommand\Header[1]{\multicolumn{1}{c}{#1}}
    \newcommand\HeaderSeconds{t,~sec}
    \centerline{%
      \begin{tabular}{lc rrr rrr}
        \toprule
        && \multicolumn{3}{c}{AVR cliques}
        & \multicolumn{3}{c}{Generated cliques}
        \\
        \cmidrule(lr){3-5} \cmidrule(l){6-8}
        &&& time in \emph{s} until & time in \emph{s} to surpass
         && time in \emph{s} until & time in \emph{s} to surpass
        \\
        \Header{Dataset}
        & \Header{Instance}
        & \Header{Objective}
        & \Header{LP-Gap 0.1\,\%}
        & \Header{METAMIS}
        & \Header{Objective}
        & \Header{LP-Gap 0.1\,\%}
        & \Header{METAMIS}
        \\
        \midrule
        \multirow{17}{*}{\rot{AVR-LARGE}}
         & CR-S-L-1 & \textbf{5,653,748} & 309 & \textbf{63} & 5,646,433 & \textbf{260} & 73 \\
         & CR-S-L-2 & \textbf{5,737,640} & 289 & 198 & 5,736,183 & \textbf{259} & \textbf{174} \\
         & CR-S-L-4 & \textbf{5,734,265} & 297 & \textbf{113} & 5,731,512 & \textbf{265} & 127 \\
         & CR-S-L-6 & 3,909,318 & 110 & 34 & \textbf{3,915,519} & \textbf{94} & \textbf{28} \\
         & CR-S-L-7 & 2,011,084 & 19 & 15 & \textbf{2,011,518} & \textbf{11} & \textbf{12} \\
         & CR-T-C-1 & \textbf{4,723,108} & 50 & 13 & 4,717,956 & \textbf{39} & 13 \\
         & CR-T-C-2 & \textbf{4,945,899} & 122 & 22 & 4,942,913 & \textbf{109} & 22 \\
         & CR-T-D-4 & \textbf{4,881,646} & \textbf{170} & 34 & 4,881,499 & 180 & \textbf{29} \\
         & CR-T-D-6 & \textbf{3,013,249} & 75 & 11 & 3,012,364 & \textbf{71} & \textbf{8} \\
         & CR-T-D-7 & \textbf{1,455,976} & 10 & 3 & 1,454,429 & \textbf{9} & 3 \\
         & CW-S-L-1 & \textbf{1,648,272} & 190 & \textbf{85} & 1,647,937 & \textbf{136} & 95 \\
         & CW-S-L-2 & 1,722,290 & 186 & 148 & \textbf{1,723,454} & \textbf{160} & \textbf{109} \\
         & CW-S-L-4 & \textbf{1,733,085} & 133 & \textbf{463} & 1,732,816 & \textbf{101} & 550 \\
         & CW-S-L-6 & 1,164,850 & 35 & 137 & \textbf{1,165,691} & \textbf{30} & \textbf{102} \\
         & CW-S-L-7 & 589,334 & 10 & \textbf{90} & \textbf{590,005} & \textbf{7} & 263 \\
         & CW-T-C-1 & \textbf{1,330,362} & 15 & \textbf{39} & 1,326,021 & \textbf{13} & 134 \\
         & CW-T-C-2 & \textbf{937,013} & 14 & 206 & 936,741 & \textbf{8} & \textbf{130} \\
        \midrule
        \multirow{16}{*}{\rot{AVR-MEDIUM}}
         & CW-T-D-4 & \textbf{460,198} & 2 & 6 & 460,010 & 2 & \textbf{3} \\
         & CW-T-D-6 & 459,776 & 2 & \textbf{7} & \textbf{459,909} & 2 & 29 \\
         & MR-D-03 & \textbf{1,758,074,154} & 1 & $<$ 1 & 1,757,878,476 & $<$ 1 & $<$ 1 \\
         & MR-D-05 & \textbf{1,790,179,623} & 1 & 1 & 1,790,072,991 & $<$ 1 & $<$ 1 \\
         & MR-D-FN & 1,802,268,039 & 1 & 1 & \textbf{1,802,369,969} & $<$ 1 & $<$ 1 \\
         & MR-W-FN & \textbf{5,386,842,781} & $<$ 1 & $<$ 1 & \textbf{5,386,842,781} & $<$ 1 & $<$ 1 \\
         & MT-D-200 & \textbf{287,155,108} & $<$ 1 & $<$ 1 & \textbf{287,155,108} & $<$ 1 & $<$ 1 \\
         & MT-D-FN & \textit{\textbf{290,866,943}} & $<$ 1 & 165 & 290,800,404 & $<$ 1 & --- \\
         & MT-W-200 & \textbf{384,056,012} & $<$ 1 & $<$ 1 & \textbf{384,056,012} & $<$ 1 & $<$ 1 \\
         & MT-W-FN & \textbf{390,869,891} & $<$ 1 & $<$ 1 & \textbf{390,869,891} & $<$ 1 & $<$ 1 \\
         & MW-D-20 & 526,688,093 & $<$ 1 & 1 & \textbf{526,931,511} & $<$ 1 & 1 \\
         & MW-D-40 & \textbf{538,944,498} & $<$ 1 & 2 & 538,769,213 & $<$ 1 & 2 \\
         & MW-D-FN & \textbf{545,720,496} & $<$ 1 & $<$ 1 & 545,432,177 & $<$ 1 & 2 \\
         & MW-W-05 & \textbf{1,329,653,408} & \textbf{18} & \textbf{1,699} & 1,329,638,653 & 19 & 1,907 \\
         & MW-W-10 & \textbf{1,300,124,841} & \textbf{18} & --- & 1,278,217,690 & 21 & --- \\
         & MW-W-FN & 1,248,553,336 & \textbf{17} & --- & \textbf{1,248,966,779} & 24 & --- \\
        \midrule
        \multirow{5}{*}{\rot{AVR-SMALL}}
         & MR-D-01 & \textbf{1,691,358,724} & $<$ 1 & 1 & 1,691,267,110 & $<$ 1 & $<$ 1 \\
         & MT-D-01 & \textit{\textbf{238,166,485}} & $<$ 1 & $<$ 1 & \textit{\textbf{238,166,485}} & $<$ 1 & $<$ 1 \\
         & MT-W-01 & \textit{\textbf{312,121,568}} & $<$ 1 & $<$ 1 & \textit{\textbf{312,121,568}} & $<$ 1 & $<$ 1 \\
         & MW-D-01 & \textbf{476,361,236} & $<$ 1 & 141 & 476,301,619 & $<$ 1 & --- \\
         & MW-W-01 & \textit{\textbf{1,270,305,952}} & $<$ 1 & $<$ 1 & \textit{\textbf{1,270,305,952}} & $<$ 1 & $<$ 1 \\
        \midrule
        \multirow{2}{*}{\rot{ALL}}
         & & \multicolumn{3}{l}{better or equal objective: \textbf{36}} & \multicolumn{3}{l}{better or equal objective: 34} \\
         & & \multicolumn{3}{l}{worse objective:           \textbf{2}}  & \multicolumn{3}{l}{          worse objective:  4}  \\
        \bottomrule
      \end{tabular}}
  \caption{
    \textbf{Impact of the clique cover on \emph{Our} method for the \emph{AVR} datasets.}
   Comparison of \cref{alg:exp-domain-stabilized} (with our primal heuristic, duality gap scheduling, and heuristic truncation, denoted \emph{Our}) applied to two different clique covers of the same \emph{AVR} instances: the cliques \emph{provided with the dataset} (\emph{AVR cliques}) and those \emph{generated greedily} from the problem graph (\emph{Generated cliques}, see \cref{sec:clique-cover-family}).
   For each construction we report the attained objective, the time until the 0.1\% relative duality gap of the respective LP relaxation is reached, and the time \emph{Our} method needs to surpass the best ``Cold''-start result of METAMIS~\cite{dong2025metaheuristic} (\cf \cref{tab:metamis-comparison}).
   For each metric, the better of the two clique covers is set in \textbf{bold} (the larger objective and the shorter runtime, respectively -- when the two coincide, both are set in bold).
   \emph{\textbf{Bold italic}} objective values coincide with those attained by METAMIS.
   All other settings are the same as in \cref{tab:metamis-comparison}.
  }
  \label{tab:metamis-clique-comparison}
\end{table}

\end{document}